\documentclass{amsart}
\title[A local limit law for anticommutators]{A local limit law for the empirical spectral distribution
of the anticommutator of independent Wigner matrices}

\author{Greg W. Anderson}
\date{August 28, 2013}
\address{School of Mathematics, University of Minnesota, Minneapolis, MN 55455 USA}
\email{gwanders@umn.edu}
\usepackage{graphicx}

\newcommand{\RRR}{{\mathcal{R}}}

\newcommand{\Ffrak}{{\mathfrak{F}}}
\newcommand{\Mfrak}{{\mathfrak{M}}}

\newcommand{\Kfrak}{{\mathfrak{K}}}
\newcommand{\Gfrak}{{\mathfrak{G}}}
\newcommand{\XXX}{{\mathcal{X}}}

\newcommand{\Ball}{{\mathrm{Ball}}}

\newcommand{\Kbold}{{\mathbf{K}}}
\newcommand{\ubold}{{\mathbf{u}}}
\newcommand{\vbold}{{\mathbf{v}}}

\newcommand{\VVV}{{\mathcal{V}}}

\newcommand{\hhh}{{\mathfrak{h}}}

\newcommand{\ii}{{\mathrm{i}}}

\newcommand{\AAA}{{\mathcal{A}}}
\newcommand{\Ibold}{{\mathbf{I}}}
\newcommand{\ebold}{{\mathbf{e}}}

\newcommand{\norm}[1]{{\left\Vert #1\right\Vert}}
\newcommand{\normnc}[1]{{\left[\!\left[ #1\right]\!\right]}}

\newcommand{\one}{{\mathbf{1}}}
\newcommand{\trace}{{\mathrm{tr}\,}}
\newcommand{\Ebold}{{\mathbf{E}}}
\newcommand{\CC}{{\mathbb{C}}}
\newcommand{\RR}{{\mathbb{R}}}

\newcommand{\Mat}{{\mathrm{Mat}}}
\newcommand{\SSS}{{\mathcal{S}}}
\newtheorem{Proposition}[subsection]{Proposition}
\newtheorem{Corollary}[subsection]{Corollary}

\newtheorem{Theorem}[subsection]{Theorem}

\newcommand{\Efrak}{{\mathfrak{E}}}

\begin{document}
\begin{abstract}
Our main result is a local limit law for  the empirical spectral distribution
of  the anticommutator of independent Wigner matrices,
modeled on the local semicircle law. 
Our approach is to adapt some techniques  from one of 
the recent papers
of Erd\"{o}s-Yau-Yin.   We also use
an algebraic description of the law of the anticommutator of free semicircular variables
due to Nica-Speicher, a self-adjointness-preserving variant of the linearization trick due to Haagerup-Schultz-Thorbj\o rnsen and the Schwinger-Dyson equation.
 A byproduct of our work is a relatively simple deterministic version of the local semicircle law.
 \end{abstract}
\maketitle
\tableofcontents

\section{Introduction and formulation of the main result}

Our main result is a local limit law for the anticommutator of independent Wigner matrices, modeled on the local semicircle law.  The latter has emerged from the recent great progress in universality for Wigner matrices.
Concerning universality, without attempting to be comprehensive, we mention
\cite{EKYY2}, \cite{EKYY}, \cite{EYY2}, \cite{EYY}, \cite{TaoVuEdge}, \cite{TaoVu} and \cite{VuWang}. The paper \cite{EYY}
has especially influenced us. We obtain our results by combining simplified variants of a few techniques from  \cite{EYY}
with variants of techniques from \cite{HST} and \cite{HT}, especially the linearization trick.
The self-adjointness-preserving variant of the linearization trick  used here
was introduced in \cite{Anderson}. (See also \cite{AndersonSlides}
and \cite{SpeicherEtAl} for slicker treatments.) It turns out to mesh well with  ``self-improving''
estimates of the type characteristic of the paper \cite{EYY}. 

\subsection{Setup for the main result}\label{section:MainResult}
We formulate our main result forthwith. See \S\ref{section:Notation} below for notation.

\subsubsection{Constants}
Fix constants $\alpha_0>0$ and $\alpha_1\geq 1$.
We also employ the absolute constant $c_{\ref{Proposition:ACNondegeneracy}}\geq 1$ defined in Proposition \ref{Proposition:ACNondegeneracy}, which is related to some special solutions of the Schwinger-Dyson equation.

\subsubsection{Random matrices}
Let $N\geq 2$ be a integer.
Let $U,V\in \Mat_N$ be random hermitian matrices with the following properties:
\begin{eqnarray}
\label{equation:Wig4}
&&\sup_{p\in [2,\infty)}p^{-\alpha_0}\left(\bigvee_{i,j=1}^N \norm{U(i,j)}_p\vee
\bigvee_{i,j=1}^N\norm{V(i,j)}_p\right)\leq \sqrt{\frac{\alpha_1}{N}}.\\
\label{equation:Wig1}
&&\mbox{The family $\{U(i,j),V(i,j)\}_{1\leq i\leq j\leq N}$ is independent.}\\
\label{equation:Wig2}
&&\mbox{All entries of $U$ and $V$ have mean zero.}\\
\label{equation:Wig3}
&&\norm{U(i,j)}_2=\norm{V(i,j)}_2=\frac{1}{\sqrt{N}}\;\;\mbox{for distinct $i,j=1,\dots,N$.}
\end{eqnarray}
This is a class of Wigner matrices  similar to that considered in \cite{EYY}.
Condition \eqref{equation:Wig4} is merely a technically convenient way of
imposing  uniformly a tail bound of exponential type. 
(See Proposition \ref{Proposition:EasyDeavy} below for the equivalence.)

\subsubsection{Apparatus from free probability} (For background see \cite[Chap. 5]{AGZ}, \cite{NicaSpeicherBis},
\cite{VDN}.)
Let $\ubold$ and $\vbold$ be freely independent semicircular noncommutative random variables.
Let $\mu_{\{\ubold\vbold\}}$ denote the law of  $\{\ubold\vbold\}$ and
let 
\begin{equation}\label{equation:macDef}
m_{\{\ubold\vbold\}}(z)=\int \frac{\mu_{\{\ubold\vbold\}}(dt)}{t-z}\;\;\;\mbox{for $z\in \hhh$}
\end{equation}
denote the Stieltjes transform of that law. Context permitting (most of the time) we will write briefly $m=m_{\{\ubold\vbold\}}(z)$.
Although $m$ depends on $z$ the notation does not show it.
It was shown in \cite[Eq. (1.15)]{NicaSpeicher} as part of a general discussion of commutators of free random variables
that $m$ satisfies the equation
\begin{equation}\label{equation:AnticommEquation}
zm^3-m^2-zm-1=0.
\end{equation}
(Caution: Our sign convention for the Stieltjes transform is opposed to that of \cite{NicaSpeicher}.)
From \eqref{equation:AnticommEquation} it follows that the support of $\mu_{\{\ubold\vbold\}}$ is $[-\zeta,\zeta]$ where
\begin{equation}\label{equation:IteratedSurd}
\zeta=\sqrt{\frac{11+5\sqrt{5}}{2}}\stackrel{\sim}{=}3.33.
\end{equation}
More precisely, it was shown that $\mu_{\{\ubold\vbold\}}$ has a density with respect to Lebesgue measure
and this density was calculated explicitly. (See  \cite[Eq. (1.17)]{NicaSpeicher}.) The density will not be needed here.

See \cite{Tetilla} for a recent discussion and application of the law  $\mu_{\{\ubold\vbold\}}$ in another context.
\subsubsection{The function $h$}
For $z\in \hhh$ let
\begin{equation}\label{equation:hdef}
h=|z+\zeta|\wedge |z-\zeta|\wedge 1.
\end{equation}
The number $0<h\leq 1$ depends on $z$ but the notation does not show it.

Here is our main result. 
\begin{Theorem}\label{Theorem:MainResult}
Notation and assumptions are as above. (Also see \S\ref{section:Notation} for general notation.) There exists a random variable $\Kbold\geq 1$
with the following two properties.
\begin{eqnarray}\label{equation:FirstKboldProperty}
&&\mbox{On the event $[\normnc{U}\vee \normnc{V}\leq 4]$ one has}\\
\nonumber&&\bigvee_{i=1}^N \left|\left(\{UV\}-z\Ibold_N\right)^{-1}(i,i)-m\right|\leq \frac{\Kbold}{\sqrt{Nh\Im z}}\\
\nonumber&&\mbox{for $z\in \hhh$ such that $|\Re z|\vee \Im z\leq 8$ and $\displaystyle\frac{4c^2_{\ref{Proposition:ACNondegeneracy}}\Kbold^2}{N}\leq h^2 \Im z$}.\\
\label{equation:SecondKboldProperty}
&&\mbox{For every $t>0$ one has $\Pr(\Kbold>t^{2\alpha_0+1})\leq \beta_0N^{\beta_1}\exp(-\beta_2t)$,}\\
\nonumber&&\mbox{for positive constants $\beta_0$ and $\beta_2$ depending only on $\alpha_0$ and $\alpha_1$}\\
\nonumber &&\mbox{and a positive  absolute constant $\beta_1$.}
\end{eqnarray}
\end{Theorem}
\noindent  The theorem is not so sharp
as the sharpest available concerning the local semicircle law.
The novelty here, rather, is to have made inroads on the general problem of proving local limit laws for  polynomials in  Wigner matrices. 
Looking forward, we have given some of our arguments in a general setting when this could be done
without making the paper significantly longer. (See \S\ref{section:Stability} and \S\ref{section:SelfConsistent} below.)
But some arguments are quite {\em ad hoc} (see \S\ref{section:Special} below)
and implicitly pose the problem of finding conceptual  general arguments with which to replace them.

One has delocalization of eigenvectors in our setup in the following  sense.
\begin{Corollary}\label{Corollary:Delocalization} Evaluate $\{UV\}$
and $\Kbold$ at a sample point of the event \linebreak $[\normnc{U}\vee \normnc{V}\leq 4]$.
We still write $\{UV\}$ and $\Kbold$ for these evaluations, respectively.
Let $\lambda$ be an eigenvalue of $\{UV\}$ and let $v$ be a corresponding unit-length (right) eigenvector.
Assume that $|\lambda|\leq 8$. Let $\rho=4c^2_{\ref{Proposition:ACNondegeneracy}}\Kbold^2/N$
and for simplicity assume that $\rho<1$.
Let $\sigma\in [\rho,\rho^{1/3}]$ be defined by the equation 
$\rho=h^2\Im z\vert_{z=\lambda+\ii\sigma}$.
Then we have
\begin{equation}\label{equation:sigmaBound}
\bigvee_{i=1}^N |v(i)|\leq \sqrt{2\sigma}.
\end{equation}
\end{Corollary}
\noindent
 This result is roughly comparable to \cite[Cor. 3.2]{EYY}.
  Figure 1 shows $\sigma$ as a function of $\lambda$
for $\rho=0.2,0.02,0.002,0.0002$.   The bound \eqref{equation:sigmaBound} is not optimal near the edge of the spectrum
and it is an open problem to optimize it.

\proof Let 
$\lambda_1\geq \cdots \geq \lambda_N$
be the eigenvalues of $\{UV\}$ and let 
$v_1,\dots,v_N$
be corresponding unit-length  eigenvectors.   We have for $i=1,\dots,N$ and $z\in \hhh$ the standard formula
$$
\frac{\Im (\{UV\}-z\Ibold_N)^{-1}(i,i)}{\Im z}=\sum_{j=1}^N \frac{|v_j(i)|^2}{|z-\lambda_j|^2}
$$
which we will apply presently.
We may assume  that 
$\lambda=\lambda_{i_0}$ and $v=v_{i_0}$
for a suitable index $i_0$.  Let
$z_0=\lambda+\ii\sigma$ and $h_0=h\vert_{z=z_0}$,
noting that 
$$|\Re z_0|\vee \Im z_0\leq 8\;\;\mbox{and}
\frac{\Kbold}{\sqrt{Nh_0\Im z_0}}= \frac{\sqrt{h_0}}{2c_{\ref{Proposition:ACNondegeneracy}}}\leq 1$$
by construction of $z_0$ and our simplifying assumption that $\rho<1$.
Thus we have 
\begin{eqnarray*}
2&\geq & 1+\frac{\Kbold}{\sqrt{Nh_0\Im z_0}}\;\geq\; \Im (\{UV\}-z_0\Ibold_N)^{-1}(i,i)\\
&=&\sum_{j=1}^N \frac{\sigma|v_j(i)|^2}{(\lambda_j-\lambda_{i_0})^2+\sigma^2}\;
\geq \;\frac{|v(i)|^2}{\sigma}
\end{eqnarray*}
by Theorem \ref{Theorem:MainResult} and the uniform bound
$|m|<1$ from Proposition \ref{Proposition:PrettyGeometry} below. 
 \qed

\begin{figure}[ht]\label{figure:Edginess}
  \caption{Closest permissible approach $\sigma$ to the real axis as a function of
$\lambda$ for  $\rho=0.2,0.02,0.002,0.0002$}
  \centering
    \includegraphics[width=0.8\textwidth]{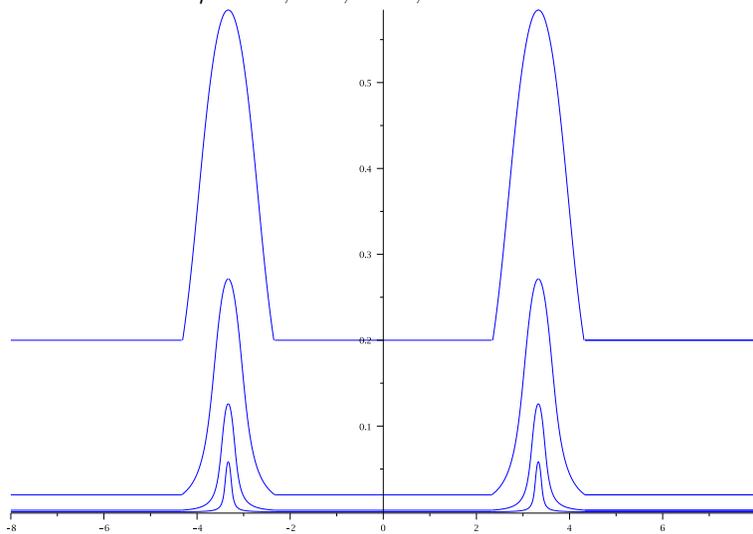}
\end{figure}

\subsection{Decay of $\Pr(\normnc{U}\vee \normnc{V}>4)$}
The conditioning on the event 
$[\normnc{U}\vee \normnc{V}\leq 4]$
taking place in Theorem \ref{Theorem:MainResult}
is not costly. In the setup of the theorem, one has
$$\Pr(\normnc{U}\vee \normnc{V}>4)\leq c_0\exp(-c_1N^{c_2})$$ for some positive constants $c_0$, $c_1$ and $c_2$
depending only on $\alpha_0$ and $\alpha_1$. 
See, e.g., the argument presented immediately after \cite[Lemma 2.1.23]{AGZ}.
The lemma in question is a combinatorial lemma somewhat weaker than the classical
result of \cite{FK} and weaker still than the more refined  results of \cite{Vu}. We will not deal further here
with the rate of decay of $\Pr(\normnc{U}\vee \normnc{V}>4)$ as $N\rightarrow\infty$.

Our proof of Theorem \ref{Theorem:MainResult} is structured overall by the following trivial remark.

\begin{Proposition}
\label{Proposition:Uroboric}
Let 
$f_1,f_2,f_3:\XXX\rightarrow[0,\infty)$ be continuous functions
on a connected topological space $\XXX$. Make the following assumptions.
\begin{eqnarray}\label{equation:Uroboric1}
&&\mbox{$f_1(x_0)<f_2(x_0)$ for some $x_0\in \XXX$.}
\\
\label{equation:Uroboric2}
&&\mbox{$f_1(x)\leq f_2(x)\Rightarrow f_1(x)\leq f_3(x)$ for all $x\in \XXX$.}\\
\label{equation:Uroboric3}
&&\mbox{$f_3(x)<f_2(x)$ for all $x\in \XXX$.}\end{eqnarray}
Then we have
\begin{equation}\label{equation:Uroboric4}
\mbox{$f_1(x)\leq f_3(x)$ for all $x\in \XXX$.}
\end{equation}
\end{Proposition}
\noindent The proposition is a less technically demanding way to think 
about estimates in the self-improving style of \cite{EYY}.
\proof We have
$\emptyset\neq \{f_1<f_2\}\subset \{f_1\leq f_3\}\subset \{f_1<f_2\}$
by hypotheses  \eqref{equation:Uroboric1}, \eqref{equation:Uroboric2} and \eqref{equation:Uroboric3}, respectively.
Since $\{f_1\leq f_3\}$ is open, closed and nonempty, in fact
$\{f_1\leq f_3\}=\XXX$ by connectedness of $\XXX$. \qed

\subsection{Further comments on methods of proofs}
\subsubsection{}
A reasonably simple explicit description of the random variable $\Kbold$ will be given later.
Given this description, the proof of property \eqref{equation:SecondKboldProperty} turns out to be an exercise
involving methods from the toolbox of \cite{EYY}. Under more restrictive hypotheses it is likely
one could obtain stronger results using the Hanson-Wright inequality. For an illuminating modern treatment of 
the latter see the recent preprint \cite{RudVer}.

\subsubsection{} The main technical result of the paper by which means we prove \eqref{equation:FirstKboldProperty}
is  a  deterministic statement of a form perhaps not seen before in connection with local limit laws. (See Theorem  \ref{Theorem:Gizmo} below.) Its proof is a reworking of the idea of a self-improving estimate---rather than marching by short steps toward the real axis, updating  estimates at each step
as in \cite{EYY}, we get our result at once by using  Proposition \ref{Proposition:Uroboric}.
\subsubsection{}We employ here  generalized resolvent techniques from \cite{Anderson}. But we do so with significant simplifications, e.g., we do not use two-variable generalized resolvents and Stieltjes transforms---rather,
we get by with just the classical parameter $z$. 
\subsection{The deterministic local semicircle law}
To facilitate comparison of our results to the literature on the local semicircle law,
we include an appendix in which we state 
and sketch a proof of the semicircular analogue of Theorem \ref{Theorem:Gizmo}.
This we call the {\em deterministic local semicircle law}.
(See Theorem \ref{Theorem:GizmoSC} below.) The latter  may be of independent interest if only for its heuristic and pedagogical value.
\subsection{Outline of the paper}
In \S\ref{section:Notation} we set out basic notation.
In \S\ref{section:QuickSD} we review the definition of the general Schwinger-Dyson equation,
including the key notion of nondegeneracy.
In \S\ref{section:Special}, which is an {\em ad hoc} mixture of free probability and  high school algebra (mostly the latter), we construct 
and analyze the particular solutions
of the Schwinger-Dyson equation needed for study of anticommutators. We also pose a problem in this section for the free probability theorists.   In \S\ref{section:Stability} we present a general stability analysis of the Schwinger-Dyson equation.
In \S\ref{section:SelfConsistent} we expeditiously analyze a matrix-valued version of the self-consistent equation \cite[Lemma 4.3]{EYY}
by algebraic and deterministic methods. (See Proposition \ref{Proposition:Frak} below.)
In \S\ref{section:GeneralizedResolvent} we do the main work
of proving \eqref{equation:FirstKboldProperty}. 
We keep the self-improving spirit of the analysis in \cite{EYY}, and continue in particular to rely heavily on (analogues of) the  formula
\begin{equation}\label{equation:MagicResolvent}
\frac{\Im(X-z\Ibold_N)^{-1}(i,i)}{\Im z}=\sum_{j=1}^N  |(X-z\Ibold_N)^{-1}(i,j)|^2
=\sum_{j=1}^N  |(X-z\Ibold_N)^{-1}(j,i)|^2
\end{equation}
where $X$ is an arbitrary $N$-by-$N$ hermitian matrix,
but our approach is deterministic and algebraic. In \S\ref{section:ProofOfMainResult}
we finish the proof of Theorem \ref{Theorem:MainResult}.
Finally, in \S\ref{section:Appendix} we briefly present our deterministic version of the local semicircle law.

\section{Notation}\label{section:Notation}
\subsection{Basic notation} Let $\{xy\}=xy+yx$ denote the anticommutator of $x$ and $y$.
We write $\ii=\sqrt{-1}$ (roman typeface).
For real numbers $x$ and $y$, let $x\vee y$ (resp., $x\wedge y$)
denote the maximum (resp., minimum) of $x$ and $y$. For $x\geq 0$, let $x_\bullet=x\vee 1$.
 Let $\Re z$ and $\Im z$ denote the real and imaginary parts of a complex number $z$, respectively,
and let $z^*$ denote the complex conjugate of $z$. 
 Let $\hhh=\{z\in \CC\mid\Im z>0\}$ denote the upper half-plane.
For a $\CC$-valued random variable $Z$ and $p\in [1,\infty)$,
let $\norm{Z}_p=(\Ebold |Z|^p)^{1/p}$ and furthermore, let $\norm{Z}_\infty$ denote the essential supremum of $|Z|$.

\subsection{Matrix notation}
Let $\Mat_{k\times \ell}$ denote the space of $k$-by-$\ell$ matrices with entries in $\CC$.
Let $\Mat_N=\Mat_{N\times N}$.
Let $\Ibold_N\in \Mat_N$ denote the  $N$-by-$N$ identity matrix. Context permitting,
we may write $1$ instead of $\Ibold_N$.
Given $A\in \Mat_{k\times \ell}$,  let $\normnc{A}$ denote the largest singular value of $A$
and let $A^*\in \Mat_{\ell\times k}$ denote the transpose conjugate  of $A$.
For $A\in \Mat_N$, let $\Re A=\frac{A+A^*}{2}$ and $\Im A=\frac{A-A^*}{2\ii}$. 
For $A\in \Mat_N$, we write $A>0$ (resp., $A\geq 0$) if $A$ is hermitian 
and positive definite (resp.,  positive semidefinite).
Given for $\nu=1,2$ a matrix $A^{(\nu)}\in \Mat_{k_\nu\times \ell_\nu}$,
recall that the {\em Kronecker product}
$A^{(1)}\otimes A^{(2)}\in \Mat_{k_1k_2\times \ell_1\ell_2}$ is defined by the rule
$$A^{(1)}\otimes A^{(2)}=\left[\begin{array}{ccc}
&\vdots&\\
\dots&A^{(1)}(i,j)A^{(2)}&\dots\\
&\vdots&
\end{array}\right].$$

\subsection{The matrix norms $\normnc{\cdot}_p$}
Given a matrix $A\in \Mat_{k\times \ell}$ with singular values 
$\mu_1\geq \mu_2\geq \cdots$
and  $p\in [1,\infty)$,
let
$\normnc{A}_p=\left(\sum_i\mu_i^p\right)^{1/p}$.
Also let  $\normnc{A}=\normnc{A}_\infty$.
Standard properties of the matrix norms
$\normnc{\cdot}_p$ are taken for granted, e.g.,  $\normnc{A}_2^2=\sum_{i,j}|A(i,j)|^2=\trace AA^*$.
Of particular importance is the
{\em H\"{o}lder inequality} which asserts that
$\normnc{AB}_r\leq \normnc{A}_p\normnc{B}_q$ whenever 
$\frac{1}{r}\leq \frac{1}{p}+\frac{1}{q}$
and the matrix product $AB$ is defined. See \cite{HoJo} or \cite{Simon} for background.
Actually only $p=1,2,\infty$ will be important.


\subsection{Stieltjes transforms}
In general, given a probability measure $\mu$ on the real line, we define its {\em Stieltjes transform} by the formula 
$S_\mu(z)=\int \frac{\mu(dt)}{t-z}$ for $z\in \hhh$.
  Note that with this sign convention we have
$\Im S_\mu(z)>0$ for $\Im z>0$. We also have a uniform bound $|S(z)|\leq 1/\Im z$.

\subsection{Inexplicit constants}
These may be denoted by $c$, $C$,  etc. and their values may change from context to context
and even from line to line. When recalling a previously defined constant we sometimes do so by
referencing as a subscript the theorem, proposition, corollary, or lemma in which the constant was defined,
e.g., $c_{\ref{Proposition:ACNondegeneracy}}$ denotes the constant
 $c$ from Proposition \ref{Proposition:ACNondegeneracy}.

\subsection{Banach spaces}\label{subsection:Banach}
Banach spaces always have complex scalars.
The norm in a Banach space $\VVV$ is denoted by $\normnc{\cdot}_\VVV$ or  simply by
$\normnc{\cdot}$ when context permits. A {\em unital} Banach algebra $\AAA$ is one equipped
with a unit $1_\AAA$ satisfying $\normnc{1_\AAA}=1$. Other notation may be used for the unit, e.g., $\Ibold_n=1_{\Mat_n}$ or $1=1_\AAA$.
 Given  a point $v_0\in \VVV$ and a constant $\epsilon\geq 0$, let  
$\Ball_\VVV(v_0,\epsilon)=\{v\in \VVV\mid \normnc{v-v_0}_\VVV\leq\epsilon\}$ (a closed ball).
Furthermore, let $B(\VVV)$ denote the space of bounded linear maps from $\VVV$ to itself
normed by the rule $\normnc{T}_{B(\VVV)}=\sup_{v\in \Ball_\VVV(0,1)}\normnc{T(v)}_\VVV$.
We invariably equip $\Mat_n$ with unital Banach algebra structure by means of
the largest-singular-value norm.



\section{A quick overview of the Schwinger-Dyson equation}\label{section:QuickSD}
For background see e.g. \cite{Anderson}, \cite{AndersonSlides}, \cite[Chap. 5]{AGZ},
\cite{HRS07} or \cite{NicaSpeicherBis}.

\subsection{Definitions}
Let $\SSS$ be a finite-dimensional unital Banach algebra.
A triple
$$(\Lambda,M,\Phi)\in \SSS\times \SSS\times B(\SSS)$$
is said to satisfy the {\em Schwinger-Dyson equation} if
\begin{equation}\label{equation:SinglePointSD}
1_\SSS+(\Lambda+\Phi(M))M=0,
\end{equation}
in which case $M$ is automatically invertible. (In a finite-dimensional unital algebra
existence of a left inverse implies existence of a two-sided inverse.) 
We emphasize that in our (slightly eccentric) usage, a solution of the Schwinger-Dyson equation is not a function;
rather, it is just a point in the space $\SSS\times \SSS\times B(\SSS)$. 
Now
let $(\Lambda,M,\Phi)$ be any solution of the  Schwinger-Dyson equation.
If the linear map
\begin{equation}\label{equation:LinearMap}
\left(x\mapsto M^{-1}x-\Phi(x)M\right)\in B(\SSS)
\end{equation}
is invertible we say that $(\Lambda,M,\Phi)$ is {\em nondegenerate} 
in which case we let $\kappa=\kappa_{\Lambda,M,\Phi}$ denote the inverse of the linear map \eqref{equation:LinearMap} and
we also say with slight abuse of terminology that the quadruple
$$(\Lambda,M,\Phi,\kappa)\in \SSS\times \SSS\times B(\SSS)\times B(\SSS)$$ is  a {\em  nondegenerate}
solution of the Schwinger-Dyson equation.  
 If we need to emphasize the role of $\SSS$
we say that $(\Lambda,M,\Phi,\kappa)$ is a solution {\em defined over $\SSS$} but we omit the epithet when context permits.
Recall our notation $x_\bullet=1\vee x$. Finally, 
we call $$\frac{1}{8\normnc{\kappa}_\bullet\normnc{\Phi}_\bullet}$$ the {\em stability radius} of $(\Lambda,M,\Phi,\kappa)$.
The meaning of the stability radius will be explained by Theorem \ref{Theorem:Stability} below.

Here is the class of examples connected to the semicircle law.
\begin{Proposition}\label{Proposition:SemicircleReview}
For $z,m\in \hhh$ satisfying 
$$z=-m^{-1}-m\;\;\left(\mbox{equivalently:}\;\;m=\frac{1}{2\pi}\int_{-2}^2\frac{\sqrt{4-t^2}\,dt}{t-z}\right)$$
one  has 
\begin{equation}\label{equation:SCreminder}
\Im m>0\;\;\mbox{and}\;\;|m|\leq 1\wedge\frac{1}{\Im z},
\end{equation}
the quadruple 
$$(z,m,1,(m^{-1}-m)^{-1})$$ is a nondegenerate solution of
the Schwinger-Dyson equation defined over $\CC$,
and the stability radius thereof satisfies the lower bound
\begin{equation}\label{equation:SCreminderBis}
\frac{1}{8|(m^{-1}-m)^{-1}|_\bullet\cdot |1|_\bullet}=\frac{1\wedge\sqrt{z^2-4}}{8}\geq\frac{\sqrt{1\wedge |z-2|\wedge |z+2|}}{8}.
\end{equation}
\end{Proposition}
\noindent This statement needs no proof.

\section{Special nondegenerate solutions of the Schwinger-Dyson equation}\label{section:Special}

For our study of anticommutators the following more exotic examples
of nondegenerate solutions of the Schwinger-Dyson equation will be needed.
\begin{Proposition}\label{Proposition:ACNondegeneracy}
Let $\Phi\in B(\Mat_3)$ be defined by the formula
\begin{equation}\label{equation:PhiIdentity}
\Phi(A)=(e_{12}+e_{21})A(e_{12}+e_{21})+(e_{13}+e_{31})A(e_{13}+e_{31})
\end{equation}
where $\{e_{ij}\}_{i,j=1}^3$ is the standard basis for $\Mat_3$.
For $z\in \hhh$, let $m=m_{\{\ubold\vbold\}}(z)\in \hhh$
be as defined on  line \eqref{equation:macDef} above. In turn, let
\begin{equation}\label{equation:CorrespondingPair}
\Lambda=\left[\begin{array}{rrr}
z&0&0\\
0&-1&0\\
0&0&1
\end{array}\right]\in \Mat_3\;\;\mbox{and}\;\;M=\left[\begin{array}{ccc}
m&0&0\\
0&-\frac{1}{m-1}&0\\
0&0&-\frac{1}{m+1}
\end{array}\right]\in \Mat_3.
\end{equation}
Although $\Lambda$ and $M$ depend on $z$ the notation does not show it. 
The triple $(\Lambda,M,\Phi)$ thus defined is a nondegenerate solution
of the Schwinger-Dyson equation defined over $\Mat_3$. 
Let
\begin{equation}\label{equation:LambdaNoughtDef}
\Lambda^0=\lim_{z\rightarrow 0}\Lambda=\left[\begin{array}{rrr}
0&0&0\\
0&-1&0\\
0&0&1\end{array}\right]\in \Mat_3.
\end{equation}
Furthermore, we have bounds
\begin{equation}\label{equation:AnticommBounds}
\normnc{\Lambda}\leq 1+|z|,\;\;
\normnc{\Phi}\leq 8,\;\;
\normnc{M+\Lambda^0}\leq 2\wedge \frac{8}{\Im z}\;\;\mbox{and}\;\;
\normnc{M}\leq 2.
\end{equation}
Let $$\kappa=\kappa_{\Lambda,M,\Phi}\in B(\Mat_3).$$
Finally, the nondegenerate solution $(\Lambda,M,\Phi,\kappa)$ 
of the Schwinger-Dyson equation 
has stability radius satisfying the lower bound
\begin{equation}\label{equation:AnticommKappaBound}
\frac{1}{8\normnc{\kappa}_\bullet\normnc{\Phi}_\bullet}\geq \frac{\sqrt{h}}{c}
\end{equation}
where $h$ is as defined on line \eqref{equation:hdef} above and $c\geq 1$ is an absolute constant.
\end{Proposition}
\noindent 
 The proof will be given in \S\ref{subsection:ProofOfACNondegeneracy} below after we have
introduced appropriate algebraic tools.  Note that the linear map $\kappa$ depends on $z$ just as do $\Lambda$ and $M$, but the notation does not show it. 
In \S\ref{subsection:PhiRecog} we provide motivation for the choice of the linear map $\Phi\in B(\Mat_3)$.
It is admittedly a flaw of paper organization that this explanation is so long postponed.
But fortunately, only the bare statement of Proposition \ref{Proposition:ACNondegeneracy} is needed in the sequel.
Thus the reader eager to see the big picture could immediately skip ahead to the next section after reading the statement of the proposition.

\subsection{Formulation of results on the equation \eqref{equation:AnticommEquation}} 
To a large extent the proof of Proposition \ref{Proposition:ACNondegeneracy} boils down to a study of
 equation \eqref{equation:AnticommEquation}. We prepare for stating several results on that equation as follows.
Let 
$$
\omega=\sqrt{\sqrt{5}-2}\stackrel{\sim}{=}\;0.4858682712,
$$
which is the unique positive root of the polynomial
\begin{equation}\label{equation:Resultantz}
m^4+4m^2-1.
\end{equation}
Let
\begin{eqnarray}
D&=&\left\{u+\ii v\bigg\vert\begin{array}{l}
u,v\in \RR,\;|u|\leq \omega\;\mbox{and}\\
0\leq v\leq ((1-4u^2)^{1/2}-u^2)^{1/2}
\end{array}\right\}\\
\nonumber&\subset&\{w\in \CC\mid \Im w\geq 0,\;|\Re w|\leq \omega\;\mbox{and}\;|w|\leq 1\}\subset \CC
\end{eqnarray}
and let $D^o$ denote the interior of  $D$.
Repeating \eqref{equation:IteratedSurd} for the reader's convenience let
\begin{equation}\label{equation:zetadef}
\zeta=\sqrt{\frac{11+ 5 \sqrt{5}}{2}}
\stackrel{\sim}{=}3.330190676,
\end{equation}
which is the unique positive root of the polynomial
\begin{equation}\label{equation:Resultantw}
z^4-11z^2-1.
\end{equation}

The algebraic results we are going to prove are as follows.
In these results and their proofs we forget about Stieltjes transforms. Instead,
we focus on the equation \eqref{equation:AnticommEquation} above and its equivalent expression \eqref{equation:AnticommEquationBis} below.
\begin{Proposition}\label{Proposition:ImplicitObstruction}
The system of equations
\begin{equation}\label{equation:ImplicitObstruction}
\begin{array}{rcl}
zm^3-m^2-zm-1&=&0\\
\frac{\partial}{\partial m}\left(zm^3-m^2-zm-1\right)&=&0
\end{array}
\end{equation} 
has exactly four complex solutions, namely
\begin{equation}\label{equation:BadPoints}
(z,m)=(-\zeta,\omega),(\zeta,-\omega),(-\ii/\zeta,\ii/\omega),(\ii/\zeta,-\ii/\omega).
\end{equation}
\end{Proposition}
\noindent Thus, in particular, we have
$
\zeta=\frac{m^2+1}{m^3-m}\vert_{m=-\omega}.
$
The four points in $\CC^2$ found in Proposition \ref{Proposition:ImplicitObstruction} are where the Implicit Function Theorem fails to yield locally a solution $m=m(z)$ of \eqref{equation:AnticommEquation} depending analytically on $z$. 

\begin{Proposition}\label{Proposition:PrettyGeometry}
(i) For each $m\in \hhh$ one has $\Im \frac{m^2+1}{m^3-m}>0$ if and only if $m\in D^o$.
(ii) For each $z\in \hhh$ there exists unique $m\in \hhh$ such that
$z=\frac{m^2+1}{m^3-m}$.
\end{Proposition}
\noindent Thus, in particular,  the result of Nica-Speicher \cite{NicaSpeicher} taken for granted, $z,m\in \hhh$ satisfy $z=\frac{m^2+1}{m^3-m}$ if and only if $m=m_{\{\ubold\vbold\}}(z)$. This is parallel to the standard fact
that $z,m\in \hhh$ satisfy  $z=-m^{-1}-m$
if and only if $m=\frac{1}{2\pi}\int_{-2}^2\frac{\sqrt{4-t^2}\,dt}{t-z}$.
\begin{Proposition}\label{Proposition:UglyCalculation}
If $z,m\in \hhh$ satisfy $z=\frac{m^2+1}{m^3-m}$, then 
\begin{eqnarray}\label{equation:BadStieltjesBound}
|m|&\leq& 1\wedge\frac{4}{\Im z}\;\;\mbox{and}\\
\label{equation:UglyCalculation}
\frac{1}{|m^2-\omega^2|}&\leq &\frac{c}{\sqrt{1\wedge |z-\zeta|\wedge |z+\zeta|}}
\end{eqnarray}
where $c\geq 1$ is an absolute constant.
\end{Proposition}
\noindent The proofs of the three propositions take up the rest of this section after we have
made the application to the proof of Proposition \ref{Proposition:ACNondegeneracy}.
\subsection{Remarks}$\;$

\subsubsection{}
We do not absolutely need Proposition \ref{Proposition:ImplicitObstruction} for the proof of Proposition \ref{Proposition:ACNondegeneracy} but it is useful to prove it in order to understand algebro-geometrically
the iterated surds $\sqrt{\sqrt{5}-2}$ and $\sqrt{\frac{11+5\sqrt{5}}{2}}$ that dominate the discussion.

\subsubsection{} Since $m=m(z)$ turns out to be a Stieltjes transform by Proposition \ref{Proposition:PrettyGeometry},
the bound \eqref{equation:BadStieltjesBound} is no surprise
and the factor $4$ can be reduced to $1$. But the bound \eqref{equation:BadStieltjesBound} is easy to obtain ``bare-handed''
and so serves as a consistency check.

\subsubsection{}\label{subsubsection:PregnantRemark}The recent paper \cite{SkSh} has elucidated finer properties of the laws of self-adjoint polynomials in free semicircular variables.  Refinement of this theory
to yield in generality the analogue of Proposition \ref{Proposition:ACNondegeneracy}
would smooth the way for a proof of a general
local limit law for self-adjoint polynomials in Wigner matrices.
We overkill the proofs of Propositions \ref{Proposition:ImplicitObstruction},
\ref{Proposition:PrettyGeometry} and \ref{Proposition:UglyCalculation} below
in hope of providing a few clues for the general theory we would like to have.

\subsection{Proof of Proposition \ref{Proposition:ACNondegeneracy}
with Propositions \ref{Proposition:ImplicitObstruction}, \ref{Proposition:PrettyGeometry}
and \ref{Proposition:UglyCalculation} granted
}\label{subsection:ProofOfACNondegeneracy}
 \subsubsection{Proof that $(\Lambda,M,\Phi)$ solves the Schwinger-Dyson equation}
  We first remark that equation \eqref{equation:AnticommEquation} 
  can be rewritten as
 \begin{equation}\label{equation:AnticommEquationBis}
 z=\frac{m^2+1}{m^3-m}=\frac{1}{m-1}+\frac{1}{m+1}-\frac{1}{m}.
 \end{equation}
 Recall that
$$
\Lambda=\left[\begin{array}{rrr}
z&0&0\\
0&-1&0\\
0&0&1
\end{array}\right]\;\;\mbox{and}\;\;M=\left[\begin{array}{ccc}
m&0&0\\
0&-\frac{1}{m-1}&0\\
0&0&-\frac{1}{m+1}
\end{array}\right].
$$
We then have
\begin{eqnarray*}
\Phi(M)&=&\left[\begin{array}{ccc}
-\frac{1}{m-1}-\frac{1}{m+1}&0&0\\
0&m&0\\
0&0&m
\end{array}\right]\;\;\mbox{and}\\
\Lambda+\Phi(M)&=&
\left[\begin{array}{ccc}
-\frac{1}{m}&0&0\\
0&m-1&0\\
0&0&m+1
\end{array}\right]\;=-M^{-1}.
\end{eqnarray*}
Thus $(\Lambda,M,\Phi)$ is indeed a solution of the Schwinger-Dyson equation.
\subsubsection{Proof of the bounds \eqref{equation:AnticommBounds}}
The first bound is clear. The second bound is proved as follows:
$$\normnc{\Phi(A)}\leq (\normnc{e_{12}+e_{21}}^2+\normnc{e_{13}+e_{31}}^2)\normnc{A}\leq 8\normnc{A}.$$
The third bound is equivalent to
$$
|m|\vee \left|\frac{m}{m-1}\right|\vee \left|\frac{m}{m+1}\right|\leq 2\left(1\wedge \frac{4}{\Im z}\right),
$$
and the latter follows easily from Propositions \ref{Proposition:PrettyGeometry} 
and \ref{Proposition:UglyCalculation}. 
Finally the fourth bound follows directly from Proposition \ref{Proposition:PrettyGeometry}.

\subsubsection{Proof of nondegeneracy}
A straightforward calculation shows that the definition \eqref{equation:PhiIdentity} can be rewritten
\begin{equation}\label{equation:PhiSecondDef}
\Phi\left(\left[\begin{array}{ccc}
x_1&x_4&x_6\\
x_5&x_2&x_8\\
x_7&x_9&x_3\end{array}\right]\right)=
\left[\begin{array}{ccc}
x_2+x_3&x_5&x_7\\
x_4&x_1&0\\
x_6&0&x_1
\end{array}\right].
\end{equation}
Abusing notation since we haven't yet proved invertibility, let $\kappa^{-1}$ denote the linear map \eqref{equation:LinearMap}.
Then we have
\begin{eqnarray*}
&&\kappa^{-1}\left(\left[\begin{array}{ccc}
x_1&x_4&x_6\\
x_5&x_2&x_8\\
x_7&x_9&x_3
\end{array}\right]\right)\\
&=&\left[\begin{array}{ccc}
1/m&0&0\\
0&-(m-1)&0\\
0&0&-(m+1)
\end{array}\right]\left[\begin{array}{ccc}
x_1&x_4&x_6\\
x_5&x_2&x_8\\
x_7&x_9&x_3
\end{array}\right]\\
&&-\left[\begin{array}{ccc}
x_2+x_3&x_5&x_7\\
x_4&x_1&0\\
x_6&0&x_1
\end{array}\right]\left[\begin{array}{ccc}
m&0&0\\
0&-\frac{1}{m-1}&0\\
0&0&-\frac{1}{m+1}
\end{array}\right].\end{eqnarray*}
With respect to the basis for $\Mat_3$ dual to the peculiar numbering
of matrix entries in \eqref{equation:PhiSecondDef}, the matrix for $\kappa^{-1}$
is block diagonal with diagonal blocks
\begin{eqnarray}\label{equation:BlockList}
&&\left[\begin{array}{ccc}
1/m&-m&-m\\
\frac{1}{m-1}&-(m-1)&0\\
\frac{1}{m+1}&0&-(m+1)
\end{array}\right],\;\;
\left[\begin{array}{cc}
1/m&\frac{1}{m-1}\\
-m&-(m-1)
\end{array}\right],\\
\nonumber&&\left[\begin{array}{cc}
1/m&\frac{1}{m+1}\\
-m&-(m+1)
\end{array}\right],\;\;\left[\begin{array}{cc}
-(m-1)&0\\
0&-(m+1)\end{array}\right],
\end{eqnarray}
respectively.
The determinants of these blocks are
\begin{equation}\label{equation:DeterminantList}
-\frac{m^4+4m^2-1}{m(m-1)(m+1)},\;\;\frac{2m-1}{m(m-1)},\;\;-\frac{2m+1}{m(m+1)},\;\;(m-1)(m+1),
\end{equation}
respectively. By Proposition \ref{Proposition:PrettyGeometry} none of the rational functions
of $m$ on the list \eqref{equation:DeterminantList} vanishes (or has a pole) in the open set $D^o$.  Thus 
$(\Lambda,M,\Phi)$ is nondegenerate and hence $\kappa$ well-defined. 

\subsubsection{Proof of the bound \eqref{equation:AnticommKappaBound}}

The inverses of the diagonal blocks on the list \eqref{equation:BlockList} are
\begin{eqnarray*}
&&\frac{\left[\begin{array}{ccc}
-(m^2-1)^2m&m^2(m^2-1)(m+1)&m^2(m^2-1)(m-1)\\
-(m+1)^2m&(2m+1)(m-1)&m^2(m+1)\\
-(m-1)^2m&m^2(m-1)&-(2m-1)(m+1)
\end{array}\right]}{(m^4+4m^2-1)},\\\\
&&\frac{\left[\begin{array}{cc}
-(m-1)^2m&-m\\
m^2(m-1)&m-1
\end{array}\right]}{2m-1},\;\;
\frac{\left[\begin{array}{cc}
(m+1)^2m&m\\
-m^2(m+1)&-(m+1)
\end{array}\right]}{2m+1},\\
&&\left[\begin{array}{cc}
-\frac{1}{m-1}&0\\
0&-\frac{1}{m+1}\end{array}\right],
\end{eqnarray*}
respectively.  We have seen that the roots of $m^4+4m^2-1$ are $\pm\omega$ and $\pm i/\omega$.
Furthermore, by Proposition \ref{Proposition:PrettyGeometry} we have $|m|< 1$
and $|\Re m|<\omega$.
Thus the entries of the matrices above are bounded in absolute value by, say,
$$\frac{6}{(1-\frac{1}{\omega})^2|m^2-\omega^2|}\vee\frac{4}{1-2\omega}\vee \frac{1}{1-\omega}\leq \frac{2^9}{|m^2-\omega^2|}.$$
It follows by Proposition \ref{Proposition:WayTrivial} immediately below that we have a bound
$$\frac{2^9}{|m^2-\omega^2|}(9+4+4+2)\sqrt{3}\leq\frac{2^{15}}{|m^2-\omega^2|}$$
for $\normnc{\kappa}$.
Finally, the bound \eqref{equation:AnticommKappaBound} follows from \eqref{equation:UglyCalculation}. \qed

 \begin{Proposition}\label{Proposition:WayTrivial}
Let $\psi\in B(\Mat_n)$ be any linear map.  Let $\{e_{ij}\}_{i,j=1}^n$ be the standard basis of $\Mat_n$
consisting of elementary matrices.
Write 
$$\psi(e_{i_2j_2})=\sum_{i_1,j_1}\psi(i_1,j_1,i_2,j_2)e_{i_1j_1}$$
for scalars $\psi(i_1,j_1,i_2,j_2)$.
Then
$$
\normnc{\psi}\leq \sqrt{n}\sum_{i_1,j_1,i_2,j_2}|\psi(i_1,j_1,i_2,j_2)|.
$$
\end{Proposition}
\noindent We omit the routine proof.

\subsection{Nonlinear $D_8$-symmetry} We now commence
a rather leisurely proof of Propositions \ref{Proposition:ImplicitObstruction}, \ref{Proposition:PrettyGeometry}
and \ref{Proposition:UglyCalculation}. 
We start by observing that the rational map
\begin{equation} \label{equation:AnticommMap} 
m\mapsto \frac{m^2+1}{m^3-m}
\end{equation}
of the $m$-line to itself
commutes with the maps
\begin{itemize}
\item $m\mapsto m^*$ (reflection in the real axis),
\item $m\mapsto -m$ ($180^o$ rotation)  and 
\item $m\mapsto \ii/m$ (composition of $90^o$ rotation and inversion).
\end{itemize}
These three maps generate an eight-element nonabelian group of symmetries  centralizing the map \eqref{equation:AnticommMap}.
Just to have a convenient short catchphrase, we refer to this phenomenon as
{\em nonlinear $D_8$-symmetry} since the group in question is isomorphic to the $8$-element dihedral group $D_8$.  

\subsection{Proof of Proposition \ref{Proposition:ImplicitObstruction}}
 The resultant of the two polynomials figuring in the system \eqref{equation:ImplicitObstruction}
with respect to $z$ is the polynomial \eqref{equation:Resultantz}
of which the full set of roots is  $\{\pm \omega,\pm \ii/\omega\}$.
The resultant of the two polynomials in \eqref{equation:ImplicitObstruction}
with respect to $w$ is the polynomial \eqref{equation:Resultantw} multiplied by $-4z$
of which the full set of roots is $\{\pm \zeta,\pm \ii/\zeta,0\}$.  (The resultants are easy to calculate using
a computer algebra system.) This gives us $20$ possible
solutions for \eqref{equation:ImplicitObstruction}.
But clearly no solution of the system \eqref{equation:ImplicitObstruction} with $z=0$ exists,
cutting the number of possibilities down to $16$.
Since equation \eqref{equation:AnticommEquation}
is linear in $z$, for each $w\in \{\pm \omega,\pm \ii/\omega\}$ there is exactly
one $z\in \{\pm\zeta,\pm \ii/\zeta\}$ such that $(w,z)$ is a solution of \eqref{equation:ImplicitObstruction}.
Thus there are exactly four solutions of \eqref{equation:ImplicitObstruction}.
One can check directly that $(\omega,-\zeta)$ is a solution of \eqref{equation:ImplicitObstruction}
and finally one gets all four solutions, namely 
the four on line \eqref{equation:BadPoints}, by exploiting nonlinear $D_8$ symmetry.
\qed

\subsection{The quadrant-lifting diagram} 
Let $m=u+\ii v$ with $u$ and $v$ real.
Then for $m^3-m\neq 0$ we have formulas
\begin{eqnarray}
\Re \frac{m^2+1}{m^3-m}
&=&\frac{u(u^4+2u^2v^2+v^4-4v^2-1)}{|m^3-m|^2},\\
\label{equation:ImaginaryIdentity}
\Im \frac{m^2+1}{m^3-m}
&=&-\frac{v(v^4+2u^2v^2+u^4+4u^2-1)}{|m^3-m|^2}.
\end{eqnarray}
It follows that
\begin{eqnarray}\label{equation:BigOval}
&&\left\{m\in \CC\setminus \{-1,0,1\}\bigg\vert\Re \frac{m^2+1}{m^3-m}=0\right\}\cup\{-1,1\}\\
\nonumber&=&
\left\{\pm \left(\sqrt{\sqrt{1+4t^2}-t^2}+\ii t\right)\bigg\vert |t|\leq \frac{1}{\omega}\right\}
\cup\ii\RR,\\
\label{equation:LittleOval}
&&\left\{m\in \CC\setminus \{-1,0,1\}\bigg\vert\Im \frac{m^2+1}{m^3-m}=0\right\}\cup\{0\}\\
\nonumber&=&\left\{\pm \left(t+\ii\sqrt{\sqrt{1-4t^2}-t^2}\right)\bigg\vert |t|\leq \omega\right\}\cup\RR.
\end{eqnarray}
\begin{figure}[ht]\label{figure:AnticommDiagram}
  \caption{The quadrant lifting diagram for $m\mapsto \frac{m^2+1}{m^3-m}$}
  \centering
    \includegraphics[width=0.8\textwidth]{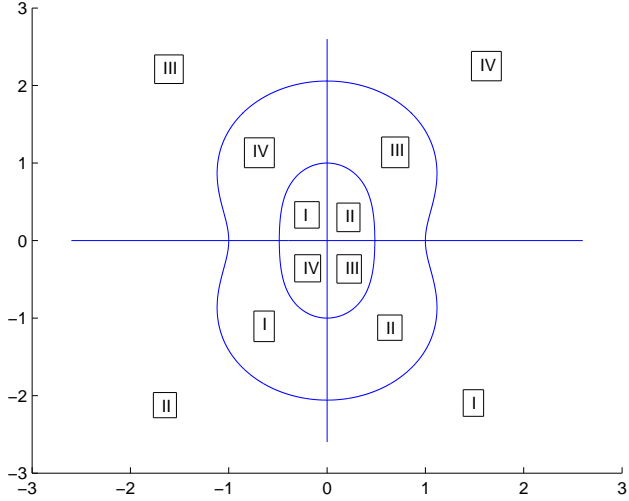}
\end{figure}
By plotting the sets \eqref{equation:BigOval} and \eqref{equation:LittleOval} on the $m$-line 
and also keeping track of the signs of $\Re \frac{m^2+1}{m^3-m}$ and $\Im \frac{m^2+1}{m^3-m}$ we obtain 
Figure 2 in which each of the twelve regions is labeled by the quadrant of the complex plane to which it is sent by the map
$m\mapsto \frac{m^2+1}{m^3-m}$.
 The diagram clearly enjoys nonlinear $D_8$-symmetry.

\subsection{Proof of Proposition \ref{Proposition:PrettyGeometry}}
The region $D$ is the union of the two regions in Figure 2 above the real axis which touch the origin.
That noted, Figure 2 (or rather, more to the point, its mode of construction)  by far overkills the proof of statement (i) of the proposition.
To prove statement (ii), in view of statement (i) already proved, it is equivalent to prove the following statement.
\begin{quote}
(iii) For each $z\in \hhh$ there exists unique $m\in D^o$ such that
$z=\frac{m^2+1}{m^3-m}$.
\end{quote}
To prove statement (iii) one begins by observing that the contour $\partial D$ bounding $D$ in the $m$-line
is sent  to the circle $\RR\cup\{\infty\}$ on the $z$-line in one-to-one and onto fashion via the map $m\mapsto \frac{m^2+1}{m^3-m}$. One then comes to the desired conclusion via the Argument Principle. \qed

\subsection{Algebraic identities}
Let $a$ and $t$ be independent (commuting) algebraic variables.  The following polynomial congruences hold
modulo $a^4+4a^2-1$ and are easy to check with a computer algebra system.
\begin{eqnarray}\label{equation:Factoid1}
&&t^4-11t^2-1\bigg\vert_{t=\frac{3a^3+13a}{2}}\equiv 0,\\\nonumber\\
\label{equation:Factoid2}
&&\left(\frac{3a^3+13a}{2}\right)\left(\frac{a^3+a}{2}\right)\equiv 1,\\
\label{equation:Factoid3}
&&(t^3-t)\pm\left(\frac{a^3+a}{2}\right)(t^2+1)\equiv \left(t\pm\frac{a^3+5a}{2}\right)(t\mp a)^2.
\end{eqnarray}
In particular, it follows that
$$
\zeta=\frac{3\omega^3+13\omega}{2}\;\;\mbox{and}\;\;\frac{1}{\zeta}=\frac{\omega^3+\omega}{2}
$$
by \eqref{equation:Factoid1} and \eqref{equation:Factoid2}, respectively.
Let
$$\rho=\frac{\omega^3+5\omega}{2}\stackrel{\sim}{=}1.272019648.$$
Now let $m,z\in \hhh$ satisfy $z=\frac{m^2+1}{m^3-m}$.
Let us substitute $(t,a)=(m,\omega)$ in \eqref{equation:Factoid3} and take the product over  the two choices of signs.
We thus obtain the identity
$$
(1-z^2/\zeta^2)(m^3-m)^2=(m^2-\rho^2)(m^2-\omega^2)^2.
$$
From the latter we immediately deduce the crucial identity
\begin{equation}\label{equation:Factoid5}
\frac{1}{|m^2-\omega^2|}=
\frac{|m^2-\rho^2|^{1/2}}{|m^2-1|}
\frac{1}{|m|}\frac{\zeta}{|z^2-\zeta^2|^{1/2}}.
\end{equation}

\subsection{Proof of Proposition \ref{Proposition:UglyCalculation}}
By Proposition \ref{Proposition:PrettyGeometry} 
we have $m\in D^o$  and hence we have crude bounds
\begin{equation}\label{equation:CrudeRho1}
|m\pm \rho|\leq 1+\rho\leq 3\;\;\mbox{and}\;\;\frac{1}{|m\pm 1|}\leq \frac{1}{1-\omega}\leq 2.
\end{equation}
Thus by rewriting \eqref{equation:ImaginaryIdentity} we obtain the relation
$$
\frac{1}{\Im \frac{m^2+1}{m^3-m}}
=\frac{|m^2-1|}{v(1-(v^4+2u^2v^2+u^4+4u^2))}|m|
\geq\frac{|m|}{4}
$$
which proves \eqref{equation:BadStieltjesBound}.
From the partial fraction expansion noted in \eqref{equation:AnticommEquationBis},
we deduce a bound
\begin{equation}\label{equation:CrudeRho2}
\frac{1}{|m|}\leq |z|+\frac{2}{1-\omega}\leq |z|+4.
\end{equation}And we have seen that $\zeta\leq 4$. 
Bounding the right side of
 \eqref{equation:Factoid5} 
by means of \eqref{equation:CrudeRho1} and \eqref{equation:CrudeRho2} we find that
\begin{equation}\label{equation:AlmostFactoid}
\frac{1}{|m^2-\omega^2|}\leq\frac{3\cdot 2^2\cdot (4+|z|)\cdot 4}{|z^2-\zeta^2|^{1/2}}.
\end{equation}
Finally, we have crude bounds
\begin{eqnarray*}
|z|\geq 10&\Rightarrow &\frac{(4+|z|)}{\sqrt{|z^2-\zeta^2|}}\leq \frac{(1+4/|z|)}{\sqrt{1-\zeta^2/|z|^2}}\leq 4,\\
|z|\geq 10&\Rightarrow &|z-\zeta|\wedge |z+\zeta|>1,\\
|z|<10&\Rightarrow &\frac{(4+|z|)}{\sqrt{|z^2-\zeta^2|}}
\leq \frac{14}{\sqrt{1\wedge |z-\zeta|\wedge |z+\zeta|}}
\end{eqnarray*}
which together imply a bound
$$\frac{4+|z|}{\sqrt{|z^2-\zeta^2|}}\leq \frac{14}{\sqrt{1\wedge |z-\zeta|\wedge |z+\zeta|}}.$$
The latter in conjunction with \eqref{equation:AlmostFactoid} proves \eqref{equation:UglyCalculation}.
\qed


\section{Stability of a general form of the Schwinger-Dyson equation}\label{section:Stability}
The main result of this section is the following. Recall our notation $x_\bullet=1\vee x$.
\begin{Theorem}[Stability of the Schwinger-Dyson equation]\label{Theorem:Stability}
Let $\SSS$ be a finite-dimensional unital Banach algebra.
Let $(\Lambda_0,M_0,\Phi_0,\kappa_0)$ be a nondegenerate solution of the Schwinger-Dyson equation
defined over $\SSS$.
Fix $G_0\in \SSS$ and let
$$E_0=1_\SSS+(\Lambda_0+\Phi_0(G_0))G_0\in \SSS.$$
We then have
\begin{equation}
\label{equation:Stability}
\normnc{G_0-M_0}\leq \frac{1}{8\normnc{\kappa_0}_\bullet\normnc{\Phi_0}_\bullet}
\Rightarrow \normnc{G_0-M_0}
\leq 20\normnc{\kappa_0}_\bullet\normnc{\Phi_0}_\bullet\normnc{M_0}_\bullet^2\normnc{E_0}.
\end{equation}
\end{Theorem}\noindent 
Statement \eqref{equation:Stability} provides the interpretation of the stability radius
$\frac{1}{8\normnc{\kappa_0}_\bullet\normnc{\Phi_0}_\bullet}$ we promised earlier to give.
The proof takes up the rest of this section after we have stated some corollaries.
The zero-subscripted notation is ugly but it helps us avoid
collisions of notation.  
Only the bare statement of the theorem is needed in the sequel
and actually only the following corollaries are  needed.

The next corollary gives the semicircular specialization of the theorem.
\begin{Corollary}\label{Corollary:SemicircularStability}
We continue in the setting of Proposition \ref{Proposition:SemicircleReview}.
Let $g,e\in \CC$
satisfy
$$e=1+(z+g)g.$$
Then we have
\begin{equation}
\label{equation:SCstability}
|g-m|\leq \frac{\sqrt{1\wedge|z-2|\wedge |z+2|}}{8}\Rightarrow
|g-m|\leq \frac{20|e|}{\sqrt{1\wedge|z-2|\wedge |z+2|}}.
\end{equation}
\end{Corollary}
\noindent This is roughly comparable to \cite[Lemma 5.2]{EYY}.
\proof Consider the instance
$$(\SSS,\Lambda_0,M_0,\Phi_0,\kappa_0,G_0,E_0)=(\CC,z,m,1,(m^{-1}-m)^{-1},g,e)$$
of Theorem \ref{Theorem:Stability}. We have
\begin{eqnarray*}
&&|g-m|\leq \frac{\sqrt{1\wedge |z+2|\wedge |z-2|}}{8}\leq \frac{1}{8|(m^{-1}-m)^{-1}|_\bullet|1|_\bullet}\\
&\Rightarrow&|g-m|\leq
20|(m^{-1}-m)^{-1}|_\bullet |1|_\bullet |e|
\leq \frac{20|e|}{\sqrt{1\wedge|z-2|\wedge |z+2|}},
\end{eqnarray*}
which proves the result.
\qed

The specialization of Theorem \ref{Theorem:Stability} relevant to the proof of Theorem \ref{Theorem:MainResult} is the following.
\begin{Corollary}\label{Corollary:StabilityAC}
We continue in the setting of Proposition \ref{Proposition:ACNondegeneracy}.
Suppose that \linebreak $G,E\in \Mat_3$ satisfy 
$$E=\Ibold_3+(\Lambda+\Phi(G))G.$$
 Then we have
 $$
\normnc{G-M}\leq \frac{\sqrt{h}}{c_{\ref{Proposition:ACNondegeneracy}}}
\Rightarrow\normnc{G-M}\leq \frac{10c_{\ref{Proposition:ACNondegeneracy}}\normnc{E}}{\sqrt{h}}.
$$
  \end{Corollary}
  \proof From Proposition \ref{Proposition:ACNondegeneracy} recall that
$$
\frac{1}{8 \normnc{\kappa}_\bullet\normnc{\Phi}_\bullet}\geq  \frac{\sqrt{h}}{c_{\ref{Proposition:ACNondegeneracy}}}\;\;
\mbox{and}\;\;\normnc{M}\leq 2,
$$
  where $c_{\ref{Proposition:ACNondegeneracy}}$ is an absolute constant.
Now consider the instance
$$(\SSS,\Lambda_0,M_0,\Phi_0,\kappa_0,G_0,E_0)=
(\Mat_3,\Lambda,M,\Phi,\kappa,G,E)
$$
of Theorem \ref{Theorem:Stability}. We have
\begin{eqnarray*}
&&\normnc{G-M}\leq \frac{\sqrt{h}}{c_{\ref{Proposition:ACNondegeneracy}}}\leq \frac{1}{8\normnc{\kappa}_\bullet\normnc{\Phi}_\bullet}\\
&\Rightarrow&\normnc{G-M}\leq
20\normnc{\kappa}_\bullet\normnc{\Phi}_\bullet\normnc{M}_\bullet^2\normnc{E}
\leq\frac{10c_{\ref{Proposition:ACNondegeneracy}}}{\sqrt{h}}\normnc{E},
\end{eqnarray*}
which proves the result. \qed

\subsection{Abbreviated terminology for the proof of Theorem \ref{Theorem:Stability}}
Until the end of the proof the linear map $\Phi_0\in B(\SSS)$ remains fixed.
Accordingly, we drop reference to $\Phi_0$ in the terminology,
saying, e.g., that the triple $(\Lambda_1,M_1,\kappa_1)$ is a nondegenerate solution of the Schwinger-Dyson equation
if the quadruple $(\Lambda_1,M_1,\Phi_0,\kappa_1)$ is. 

\subsection{The deformation equation associated to a nondegenerate solution of the Schwinger-Dyson equation}
As in the statement of Theorem \ref{Theorem:Stability}, let
$$(\Lambda_0,M_0,\kappa_0)\in \SSS\times \SSS\times B(\SSS)$$ be a nondegenerate solution of the Schwinger-Dyson equation.
We say that a pair $(\Theta,H)\in \SSS\times \SSS$
satisfies the {\em deformation equation} associated with the triple $(\Lambda_0,M_0,\kappa_0)$ if
\begin{equation}\label{equation:LocalSchwingerDyson}
H=\kappa_0\left(\Theta M_0+\Theta H+\Phi_0(H)H\right).
\end{equation}

\begin{Proposition}\label{Proposition:LocalSchwingerDyson}
As in the statement of Theorem \ref{Theorem:Stability}, let
$(\Lambda_0,M_0,\kappa_0)$ be a nondegenerate solution of the Schwinger-Dyson equation. 
Fix $(\Lambda_1,M_1)\in \SSS\times \SSS$ and
write 
$(\Theta,H)=(\Lambda_1-\Lambda_0,M_1-M_0)$.
 Then 
the pair $(\Lambda_1,M_1)$ is a solution of the Schwinger-Dyson equation
if and only if the pair $(\Theta,H)$ is a solution of the deformation equation \eqref{equation:LocalSchwingerDyson}
associated with the triple $(\Lambda_0,M_0,\kappa_0)$. \end{Proposition}
\proof We first prove the implication ($\Rightarrow$).
We have
\begin{eqnarray*}
0&=&1+(\Lambda_1+\Phi_0(M_1))M_1\;=\;1+(\Lambda_0+\Theta+\Phi_0(M_0+H))(M_0+H)\\
&=&1+(\Lambda_0+\Phi_0(M_0))M_0+(\Theta+\Phi_0(H))H\\
&&+(\Theta+\Phi_0(H))M_0+(\Lambda_0+\Phi_0(M_0))H\\
&=&\Theta H+\Theta M_0+\Phi_0(H)H+\Phi_0(H)M_0-M_0^{-1}H
\end{eqnarray*}
and hence
$$M_0^{-1}H-\Phi_0(H)M_0\;=\;\Theta M_0+\Theta H+\Phi_0(H)H.$$
Thus the deformation equation \eqref{equation:LocalSchwingerDyson} holds.  The steps of the preceding argument are reversible. Thus the converse  ($\Leftarrow$) also holds.
\qed

\begin{Proposition}\label{Proposition:SeparationDeformation}
As in the statement of Theorem \ref{Theorem:Stability},
let $(\Lambda_0,M_0,\kappa_0)$ be a nondegenerate solution of the  Schwinger-Dyson equation.
Fix constants $\epsilon$ and $\delta$ such that
$$
0\leq \epsilon \leq \frac{1}{4\normnc{\kappa_0}_\bullet\normnc{\Phi_0}_\bullet}\;\;\mbox{and}\;\;
0\leq \delta \leq \frac{\epsilon}{4\normnc{\kappa_0}_\bullet\normnc{M_0}_\bullet}.
$$
Fix $\Lambda\in \Ball_\SSS(\Lambda_0,\delta)$. (For this notation see \S\ref{subsection:Banach}.)
Then there exists unique \linebreak $M\in \Ball_\SSS(M_0,\epsilon)$ such that the pair
$(\Lambda,M)$ is a solution of the Schwinger-Dyson equation. 
\end{Proposition}
\proof
Let  $$\Theta=\Lambda-\Lambda_0\in \Ball_\SSS(0,\delta)$$
and consider the quadratic mapping
$$Q:=\left(x\mapsto\kappa_0\left(\Theta M_0+\Theta x+\Phi_0(x)x\right)\right):\SSS\rightarrow\SSS.$$
By Proposition \ref{Proposition:LocalSchwingerDyson}, 
an element $M\in \SSS$ has the property that the pair $(\Lambda,M)$ is a solution of the Schwinger-Dyson equation
if and only if the difference $M-M_0$ is a fixed point of $Q$. Thus our task is transformed to that of proving
the existence of a unique fixed point of $Q$ in $\Ball_\SSS(0,\epsilon)$. For achieving the latter goal the Banach fixed point theorem is the natural tool.

We turn now to the analysis of $Q$ restricted to $\Ball_\SSS(0,\epsilon)$.
For  $x\in \Ball_\SSS(0,\epsilon)$  we have
\begin{eqnarray*}
\normnc{Q(x)}&=&\normnc{\kappa_0\left(\Theta M_0+\Theta x+\Phi_0(x)x\right)}\\
&\leq &\normnc{\kappa_0}\normnc{M_0}\delta+\normnc{\kappa_0}\delta\epsilon+
\normnc{\kappa_0}\normnc{\Phi_0}\epsilon^2\;\leq \;\frac{\epsilon}{4}+\frac{\epsilon}{4}+\frac{\epsilon}{4}\leq \epsilon.
\end{eqnarray*}
Thus we have
\begin{equation}
\label{equation:Capture}
Q\left(\Ball_\SSS(0,\epsilon)\right)\subset \Ball_\SSS(0,\epsilon).
\end{equation}
For $x_1,x_2\in \Ball_\SSS(0,\epsilon)$ we have
\begin{eqnarray*}
&&\normnc{Q(x_1)-Q(x_2)}\\
&=&\normnc{\kappa_0\left(\Theta M_0+\Theta x_1+\Phi_0(x_1)x_1\right)
-\kappa_0\left(\Theta M_0+\Theta x_2+\Phi_0(x_2)x_2\right)}\\
&\leq &\normnc{\kappa_0}\normnc{\Theta(x_1-x_2) +\Phi_0(x_1-x_2)x_1
+\Phi_0(x_2)(x_1-x_2)}\\
&\leq &\left(\normnc{\kappa_0}\delta+\normnc{\kappa_0}\normnc{\Phi_0}\epsilon+\normnc{\kappa_0}\normnc{\Phi_0}\epsilon\right)\normnc{x_1-x_2}\\
&\leq&\left(\frac{1}{4}+\frac{1}{4}+\frac{1}{4}\right)\normnc{x_1-x_2}=\frac{3}{4}\normnc{x_1-x_2}.
\end{eqnarray*}
Thus we have
\begin{equation}
\label{equation:Contract}
x_1,x_2\in \Ball_\SSS(0,\epsilon)\Rightarrow \normnc{Q(x_1)-Q(x_2)}\leq \frac{3}{4}\normnc{x_1-x_2}.
\end{equation}
By \eqref{equation:Capture}  and \eqref{equation:Contract} 
the map $Q$  induces a contraction mapping of the complete metric space $\Ball_\SSS(0,\epsilon)$ to itself.
By the Banach fixed point theorem $Q$ indeed has a unique fixed point in $\Ball_\SSS(0,\epsilon)$.
\qed

\subsection{Proof of Theorem \ref{Theorem:Stability}}\label{subsection:ProofOfStability}
$$\normnc{G_0-M_0}\leq \frac{1}{8\normnc{\kappa_0}_\bullet\normnc{\Phi_0}_\bullet}
\Rightarrow \normnc{G_0-M_0}
\leq 20\normnc{\kappa_0}_\bullet\normnc{\Phi_0}_\bullet\normnc{M_0}_\bullet^2\normnc{E_0}.$$
We may assume that
\begin{eqnarray}
\label{equation:StabHyp1}
\normnc{E_0}&\leq&\frac{1}{64\normnc{\kappa_0}_\bullet^2\normnc{M_0}_\bullet^2\normnc{\Phi_0}_\bullet^2},\end{eqnarray}
since otherwise \eqref{equation:Stability} holds automatically and there is nothing to prove.
Now by the hypothesis of \eqref{equation:Stability} we have $\normnc{G_0}\leq2\normnc{M_0}_\bullet$
and furthermore by \eqref{equation:StabHyp1} we have 
$\normnc{E_0}\leq \frac{1}{2}$.
Thus $\Lambda_0+\Phi_0(G_0)$ is invertible and its inverse satisfies the bound
\begin{equation}\label{equation:StabCons1}
\normnc{(\Lambda_0+\Phi_0(G_0))^{-1}}\leq  2\normnc{G_0}\leq 4\normnc{M_0}_\bullet.\end{equation}
Let
$$M=-(\Lambda_0+\Phi_0(G_0))^{-1}\;\;\mbox{and}\;\;\Lambda=\Lambda_0+\Phi_0(G_0-M).$$
 The pair $(\Lambda,M)$ is a solution of the Schwinger-Dyson equation 
 because
 \begin{eqnarray*}
1+(\Lambda+\Phi_0(M))M
&=&1+(\Lambda_0+\Phi_0(G_0-M)+\Phi_0(M))M\\
&=&1+(\Lambda_0+\Phi_0(G_0))M=1-1=0.
\end{eqnarray*}
 By \eqref{equation:StabCons1} and the definitions we have
 \begin{eqnarray}\label{equation:StabCons1bis}
 \normnc{G_0-M}&=&\normnc{(\Lambda_0+\Phi_0(G_0))^{-1}+G_0}\\
\nonumber &=&\normnc{(\Lambda+\Phi_0(G_0))^{-1}E_0}\;\leq\; 4\normnc{M_0}_\bullet \normnc{E_0}.
 \end{eqnarray}
By hypothesis of \eqref{equation:Stability} along with \eqref{equation:StabHyp1}
and \eqref{equation:StabCons1bis} we have
\begin{eqnarray}\label{equation:StabCons3}
\normnc{M-M_0}&\leq&\normnc{G_0-M}+\normnc{G_0-M_0}\\
\nonumber&\leq& 4\normnc{M_0}_\bullet\normnc{E_0}+\frac{1}{8\normnc{\kappa_0}_\bullet\normnc{\Phi_0}_\bullet}\leq \frac{1}{4\normnc{\kappa_0}_\bullet\normnc{\Phi_0}_\bullet}.\end{eqnarray}
By \eqref{equation:StabHyp1} and \eqref{equation:StabCons1bis} we also have
\begin{equation}\label{equation:StabCons2}
\normnc{\Lambda-\Lambda_0}\leq  4\normnc{\Phi_0}_\bullet \normnc{M_0}_\bullet \normnc{E_0}\leq \frac{1}{16\normnc{\kappa_0}_\bullet^2\normnc{M_0}_\bullet\normnc{\Phi_0}_\bullet}.\end{equation}
Applying Proposition \ref{Proposition:SeparationDeformation} in the case 
$$(\delta,\epsilon)=\left(\frac{1}{16\normnc{\kappa_0}_\bullet^2\normnc{M_0}_\bullet\normnc{\Phi_0}_\bullet},
\frac{1}{4\normnc{\kappa_0}_\bullet\normnc{\Phi_0}_\bullet}\right),$$ we conclude that $M$ is the unique
element of $\Ball_\SSS\left(M_0,\frac{1}{4\normnc{\kappa_0}_\bullet\normnc{\Phi_0}_\bullet}\right)$ such that $(\Lambda,M)$ is a solution of
the Schwinger-Dyson equation. By applying Proposition \ref{Proposition:SeparationDeformation}
again in the case 
$$(\delta,\epsilon)=(\normnc{\Lambda-\Lambda_0},4\normnc{\kappa_0}_\bullet\normnc{M_0}_\bullet\normnc{\Lambda-\Lambda_0})$$ we find that in fact 
$$\normnc{M-M_0}\leq 4\normnc{\kappa_0}_\bullet\normnc{M_0}_\bullet\normnc{\Lambda-\Lambda_0}.$$
Thus by \eqref{equation:StabCons1bis} and \eqref{equation:StabCons2} we have
\begin{eqnarray*}
\normnc{G_0-M_0}&\leq &\normnc{G_0-M}+\normnc{M-M_0}\\
&\leq&
4\normnc{M_0}_\bullet\normnc{E_0}+(4\normnc{\kappa_0}_\bullet\normnc{M_0}_\bullet)(4\normnc{M_0}_\bullet\normnc{\Phi_0}_\bullet)\normnc{E_0},
\end{eqnarray*}
which suffices to prove \eqref{equation:Stability}. \qed


\section{A matrix-valued self-consistent equation}\label{section:SelfConsistent}
We prove a technical result similar in intent to \cite[Lemma 4.3]{EYY}
although rather different because instead of being probabilistic it is formal and algebraic.
(See Proposition \ref{Proposition:Frak}  below.) 
In any case, the object of study, namely the self-consistent equation, is essentially the same.

\subsection{Setup for the technical result}
Fix a finite-dimensional unital Banach algebra $\SSS$.
Fix a nondegenerate solution $$(\Lambda_0,M_0,\Phi_0,\kappa_0)$$
of the Schwinger-Dyson equation defined over $\SSS$
for which (recall)
 $\frac{1}{8\normnc{\kappa_0}_\bullet \normnc{\Phi_0}_\bullet}$
is by definition the stability radius. 
Fix a family
$$\{G_i,\widehat{G}_i\}_{i=1}^N$$
of elements of $\SSS$
where all the $G_i$ are invertible. 
Consider the statistic
\begin{eqnarray*}
\Efrak&=&\bigvee_{i=1}^N
\frac{\normnc{G_i^{-1}+\Lambda_0+\Phi_0(\widehat{G}_i)}}{\normnc{\widehat{G}_i}_\bullet^{1/2}}\vee\bigvee_{i=1}^N
\sqrt{\frac{\normnc{\widehat{G}_i-\frac{1}{N}\sum_{i=1}^N G_i}}{\normnc{G_i}_\bullet\normnc{G_i^{-1}}}},
\end{eqnarray*}
which is a gauge of error in this situation.
The idea to emphasize the statistic $\Efrak$ clearly derives from \cite[Lemma 4.3]{EYY}
and the related constellation of identities and estimates. 
\begin{Proposition}\label{Proposition:Frak} Notation and assumptions are as above. We have
\begin{eqnarray}
\label{equation:FrakBis}
&&\bigvee_{i=1}^N \normnc{G_i-M_0}\leq \frac{1}{8\normnc{\kappa_0}_\bullet\normnc{\Phi_0}_\bullet}\\
\nonumber&\Rightarrow &\bigvee_{i=1}^N \normnc{G_i-M_0}\leq 2^{14}(1+\normnc{M_0})^7(\normnc{\Phi_0}_\bullet\vee\normnc{\Lambda_0}_\bullet)^4\normnc{\kappa_0}_\bullet\Efrak.
\end{eqnarray}
\end{Proposition}
\noindent Note the similarity in form to  hypothesis \eqref{equation:Uroboric2}
of Proposition \ref{Proposition:Uroboric}. 
\proof 
Let
$$\Mfrak\;=\;1+\normnc{M_0}\;\;\mbox{and}\;\;\Ffrak=\normnc{\Phi_0}_\bullet\vee \normnc{\Lambda_0}_\bullet.$$
Let 
$$\Gfrak=\bigvee_{i=1}^N \normnc{G_i}_\bullet,\;\;G=\frac{1}{N}\sum_{i=1}^N G_i\;\;\mbox{and}\;\;E=1+(\Lambda_0+\Phi_0(G))G.$$
We temporarily beg the question by assuming
\begin{equation}
\label{equation:Frak}
\normnc{E}\vee \bigvee_{i=1}^N \normnc{G_i-G}
\leq 2^8\Gfrak^5\Ffrak^3\Efrak.
\end{equation}
By the hypothesis of \eqref{equation:FrakBis} we have 
$\normnc{G-M_0}\leq \frac{1}{8\normnc{\kappa_0}_\bullet\normnc{\Phi_0}_\bullet}$
and $\Gfrak\leq \Mfrak$. Thus we have
\begin{eqnarray*}
\bigvee_{i=1}^N \normnc{G_i-M_0}&\leq &\normnc{G-M_0}+\bigvee_{i=1}^N \normnc{G_i-G}
\;\leq \; 20\normnc{\kappa_0}_\bullet\normnc{\Phi_0}_\bullet\normnc{M_0}_\bullet^2\normnc{E}+
2^8\Gfrak^5\Ffrak^3\Efrak\\
&\leq & (2^5\normnc{\kappa_0}_\bullet\Ffrak\Mfrak^2+1)2^8\Mfrak^5\Ffrak^3\Efrak
\leq 2^{14}\Mfrak^7\Ffrak^4\normnc{\kappa_0}_\bullet\Efrak
\end{eqnarray*}
 by Theorem \ref{Theorem:Stability} and \eqref{equation:Frak},
 i.e., \eqref{equation:FrakBis} holds.

It remains now only to prove \eqref{equation:Frak}. (We will not need the hypothesis of \eqref{equation:FrakBis}
for that purpose.)
We may assume that
\begin{equation}\label{equation:NuffEfrak}
\Efrak^2\leq \Efrak\leq \frac{1}{2^6\Gfrak^3\Ffrak^2}\leq 1
\end{equation}
because the left side of \eqref{equation:Frak} is trivially bounded by $2^2\Gfrak^2\Ffrak$.

We first bound $\normnc{G-\widehat{G}_i}$. We calculate as follows. 
 \begin{eqnarray*}
\normnc{G_i^{-1}}&\leq &\normnc{G_i^{-1}+\Lambda_0+\Phi_0(\widehat{G}_i)}+\Ffrak+\Ffrak\normnc{\widehat{G}_i}\\
&\leq &\Efrak\normnc{\widehat{G}_i}_\bullet^{1/2}+2\Ffrak\normnc{\widehat{G}_i}_\bullet\leq4\Ffrak\normnc{\widehat{G}_i}_\bullet\\
&\leq &4\Ffrak\normnc{G}_\bullet+4\Ffrak\normnc{G-\widehat{G}_i}\;\leq \;4\Gfrak\Ffrak+4\Ffrak\Efrak^2\normnc{G_i}_\bullet\normnc{G_i^{-1}}\\
&\leq &4\Gfrak\Ffrak+4\Gfrak\Ffrak\Efrak \normnc{G_i^{-1}}.
\end{eqnarray*}
Since $4\Gfrak\Ffrak\Efrak\leq \frac{1}{2}$  by \eqref{equation:NuffEfrak} and
hence $\normnc{G_i^{-1}}\leq 8\Gfrak\Ffrak$
we have
\begin{equation}\label{equation:FirstGoal}
\normnc{G-\widehat{G}_i}\leq \Efrak^2 \normnc{G_i}_\bullet \normnc{G_i^{-1}}\leq 8\Gfrak^2\Ffrak\Efrak.
\end{equation}

We next bound $\normnc{E}$. We calculate as follows.
\begin{eqnarray*}
\normnc{G_i^{-1}+\Lambda_0+\Phi_0(G)}
&\leq&\Efrak\normnc{\widehat{G}_i}^{1/2}_\bullet+\Ffrak\normnc{G-\widehat{G}_i}\\
&\leq &\Efrak\normnc{G}^{1/2}_\bullet+\Efrak \normnc{G-\widehat{G}_i}^{1/2}+\Ffrak\normnc{G-\widehat{G}_i}\\
&\leq &\Efrak\normnc{G}_\bullet^{1/2}+\Efrak^2+\normnc{G-\widehat{G}_i}+\Ffrak\normnc{G-\widehat{G}_i}\\
&\leq &2\Gfrak^{1/2}\Efrak+2\Ffrak\normnc{G-\widehat{G}_i}\\
&\leq &2\Gfrak^{1/2}\Efrak+16\Gfrak^2\Ffrak^2\Efrak \leq 2^5\Gfrak^2\Ffrak^2\Efrak.
\end{eqnarray*}
We used the arithmetic-geometric mean inequality at the third step above
and \eqref{equation:FirstGoal} at the penultimate step.
We conclude that
\begin{equation}\label{equation:TechCons1}
\normnc{E}
\leq \bigvee_{i=1}^N \normnc{1+(\Lambda_0+\Phi_0(G))G_i}\leq 2^5\Gfrak^3\Ffrak^2\Efrak.
\end{equation}

Finally we bound $\normnc{G-G_i}$. By \eqref{equation:NuffEfrak}, the left side of \eqref{equation:TechCons1} is bounded by $\frac{1}{2}$.
Thus $\Lambda_0+\Phi_0(G)$ is invertible and we have
$$\normnc{(\Lambda_0+\Phi_0(G))^{-1}}\leq 2\normnc{G}.$$
In turn we have by \eqref{equation:TechCons1} that
$$
\normnc{(\Lambda_0+\Phi_0(G))^{-1}+G}\vee\normnc{(\Lambda_0+\Phi_0(G))^{-1}+G_i}
\nonumber\\
\leq 2^6\Gfrak^4\Ffrak^2\Efrak
$$
and hence
\begin{equation}\label{equation:TechCons2}
\bigvee_{i=1}^N \normnc{G-G_i}
\leq 2^7\Gfrak^4\Ffrak^2\Efrak.
\end{equation}
The bound \eqref{equation:Frak} follows now
from \eqref{equation:TechCons1} and \eqref{equation:TechCons2}.
The proof of Proposition \ref{Proposition:Frak} is complete. \qed


\section{The generalized resolvent for anticommutators}\label{section:GeneralizedResolvent}
In this section all considerations are  algebraic and deterministic except in \S\ref{subsection:PhiRecog}.
All constructions here proceed from a couple of arbitrarily chosen  hermitian matrices $U,V\in \Mat_N$
and a complex number $z\in \hhh$.
In \S\ref{section:ProofOfMainResult} we will take $U$ and $V$ to be the random matrices figuring in Theorem \ref{Theorem:MainResult} but in this section the randomness stays in the background.
The main result of this section specializes Proposition \ref{Proposition:Frak} and lays the groundwork for the construction of the random variable $\Kbold$ figuring in Theorem \ref{Theorem:MainResult}. (See Theorem \ref{Theorem:Gizmo} below.)

\subsection{The quadruple $(\Lambda,M,\Phi,\kappa)$}
Let $(\Lambda,M,\Phi,\kappa)$ be the nondegenerate solution of the Schwinger-Dyson equation defined in
Proposition \ref{Proposition:ACNondegeneracy}. The objects $\Lambda$, $M$ and $\kappa$
depend on $z\in \hhh$ but the notation does not show it. 

\subsection{The matrices $X$ and $W$}
Fix an integer $N\geq 2$.
Fix hermitian matrices $U,V\in \Mat_N$. These remain fixed throughout this section.
Let
\begin{eqnarray*}
X&=&\left[\begin{array}{ccc} 
0&\frac{U-V}{\sqrt{2}}&\frac{-U-V}{\sqrt{2}}\\
\frac{U-V}{\sqrt{2}}&0&0\\
\frac{-U-V}{\sqrt{2}}&0&0\end{array}\right]\in\Mat_{3N}\;
\mbox{and}\;\\
W&=&\left[\begin{array}{ccc}
\Ibold_N&0&0\\
\frac{-U+V}{\sqrt{2}}&\Ibold_N&0\\
\frac{-U-V}{\sqrt{2}}&0&\Ibold_N
\end{array}\right]\in \Mat_{3N}.
\end{eqnarray*}
Note that $X$ is hermitian.
Note that
\begin{equation}\label{equation:StupidXW}
1\leq \normnc{W}=\normnc{W^{-1}}=\normnc{W^*}=\normnc{(W^*)^{-1}}\;\;\mbox{and}\;\;
\normnc{X}\vee \normnc{W}\leq 8(\normnc{U}\vee \normnc{V}\vee 1).
\end{equation}
\subsection{Definition of the generalized resolvent $R$}
Note that we have for arbitrary $z\in \hhh$ a factorization
$$
W^*(X-\Lambda\otimes \Ibold_N)W=\left[\begin{array}{crr}
UV+VU-z\Ibold_N&0&0\\
0&\Ibold_N&0\\
0&0&-\Ibold_N
\end{array}\right].
$$
It follows that $X-\Lambda\otimes \Ibold_N$ is invertible
and that
\begin{equation}\label{equation:BasicRid}
R=\left(X-\Lambda\otimes \Ibold_N\right)^{-1}=
W\left[\begin{array}{crr}
(UV+VU-z\Ibold_N)^{-1}&0&0\\
0&\Ibold_N&0\\
0&0&-\Ibold_N
\end{array}\right]W^*.
\end{equation}
The {\em generalized resolvent} $R$ thus defined depends on $z$ but the notation does not show it.
Note that the resolvent of the anticommutator 
$\{UV\}$ appears as the \linebreak $N$-by-$N$ block in the upper left corner of $R$.
For discussion of the self-adjoint linearization trick whereby $R$ has been contrived see \cite{Anderson},
\cite{AndersonSlides} or \cite{SpeicherEtAl}.

\subsection{Specialized matrix notation} 
 Let $e_i\in \Mat_{1\times N}$ denote the $i^{th}$ row of $\Ibold_N$
 and let  $\hat{e}_i\in \Mat_{(N-1)\times N}$ denote $\Ibold_N$ with the $i^{th}$ row deleted.
Let  $\ebold_i=\Ibold_3\otimes e_i\in \Mat_{3\times 3N}$ and
 $\hat{\ebold}_i=\Ibold_3\otimes \hat{e}_i\in \Mat_{3(N-1)\times 3N}$.

\subsection{Objects associated with $R$}
For $i=1,\dots,N$ and $z\in \hhh$, let
\begin{eqnarray*}
G_i&=&\ebold_iR\ebold^*_i\in \Mat_3,\;\;
G\;=\;\frac{1}{N}\sum_{i=1}^NG_i\in \Mat_3,\\\label{equation:RiDef}
R_i&=&(\hat{\ebold}_iX\hat{\ebold}_i^*-\Lambda\otimes \Ibold_{N-1})^{-1}\in \Mat_{3(N-1)},\\
\widehat{G}_i&=&\frac{1}{N}\sum_{j=1}^N \ebold_j\hat{\ebold}_i^*R_i\hat{\ebold}_i\ebold_j^*\in \Mat_3,\\
\label{equation:QiDef}
Q_i&=&\ebold_iX\hat{\ebold}_i^*R_i\hat{\ebold}_iX\ebold_i^*-\ebold_iX\ebold_i^*-\Phi(\widehat{G}_i)\in \Mat_3,\\
\Kfrak_i&=&
 1\vee\frac{\normnc{Q_i}}{\frac{1}{\sqrt{N}}\left(1\vee \frac{\normnc{R_i}_2}{\sqrt{N}}\right)}
 \in [1,\infty)\;\;\mbox{and}\;\;
\Kfrak=\bigvee_{i=1}^N \Kfrak_i.
\end{eqnarray*}
All these objects depend on $z$ but the notation does not show it.
Furthermore the $z$-dependence is continuous.
Except for $\Kfrak_i$ and $\Kfrak$, the $z$-dependence is analytic.

\subsection{Structure of $X$ as a random matrix}\label{subsection:PhiRecog}
Suppose for the moment that $U$ and $V$ are random and as in Theorem \ref{Theorem:MainResult}
satisfy \eqref{equation:Wig4}, \eqref{equation:Wig1}, \eqref{equation:Wig2} and \eqref{equation:Wig3}. We claim that
the random matrix $X$ has the following properties.
\begin{eqnarray}
\label{equation:Wig4bis}
&&\sup_{p\in [2,\infty)}p^{-\alpha_0}\bigvee_{i,j=1}^N \norm{\normnc{\ebold_iX\ebold_j}}_p<\alpha_2\\
\nonumber&&\mbox{for a constant $\alpha_2$ depending only on $\alpha_0$ and $\alpha_1$.}\\
\label{equation:Wig1bis}
&&\mbox{The family $\{\ebold_iX\ebold_j\}_{1\leq i\leq j\leq N}$ is independent.}\\
\label{equation:Wig2bis}
&&\Ebold X=0.\\
\label{equation:Wig3bis}
&&\Ebold \ebold_i X\ebold_j^* A\ebold_k^*X\ebold_i=\delta_{jk}\Phi(A)\\
\nonumber&&\;\;\mbox{for $i,j,k=1,\dots,N$ s.t. $i\not\in \{j,k\}$ and $A\in \Mat_3$.}
\end{eqnarray}
The first three claims are clear. We just prove the last.
We have in any case
$$
\ebold_iX\ebold_j^*=\left(\frac{U-V}{\sqrt{2}}(i,j)\right)(e_{12}+e_{21})+\left(\frac{-U-V}{\sqrt{2}}(i,j)\right)(e_{13}+e_{31})
$$
by direct appeal to the definitions.
Now by assumptions
\eqref{equation:Wig1}, \eqref{equation:Wig2} and \eqref{equation:Wig3},
for any fixed distinct indices $i,j=1,\dots,N$, the two $\CC$-valued random variables
$$\frac{U-V}{\sqrt{2}}(i,j)\;\;\mbox{and}\;\;\frac{-U-V}{\sqrt{2}}(i,j)$$
form an orthonormal system. Formula \eqref{equation:Wig3bis} then follows by
the definition of $\Phi$.
The claims are proved. From the claims it follows  that for $i=1,\dots,N$ we have
\begin{eqnarray}\label{equation:PhiRecog}
&&\mbox{$\sigma(\hat{\ebold}_iX\hat{\ebold}_i^*)$ and $\sigma(\ebold_iX)$ are independent},\\
\label{equation:Wig6bis}
&&\mbox{$R_i$ and $\widehat{G}_i$ are $\sigma(\hat{\ebold}_iX\hat{\ebold}_i^*)$-measurable and}\\
\label{equation:Wig5bis}
&&\Ebold (Q_i\vert \hat{\ebold}_iX\hat{\ebold}_i^*)=0\;\mbox{a.s..}
\end{eqnarray}
Achievement of the property \eqref{equation:Wig5bis}
is the principal motivation for the definition of $\Phi$.
Our probabilistic digression is now concluded. We return to an algebraic viewpoint
for the rest of \S\ref{section:GeneralizedResolvent}. 
\subsection{The two-by-two  inversion formula}
We quickly review some standard algebraic gadgetry.
For a matrix $\left[\begin{array}{cc}
a&b\\
c&d\end{array}\right]\in \Mat_N$ decomposed into blocks with $a$ and $d$ square we have a factorization
$$\left[\begin{array}{cc}
a&b\\
c&d\end{array}\right]
=\left[\begin{array}{cc}
1&bd^{-1}\\
0&1
\end{array}\right]\left[\begin{array}{cc}
a-bd^{-1}c&0\\
0&d
\end{array}\right]\left[\begin{array}{cc}
1&0\\
d^{-1}c&1
\end{array}\right]
$$
whenever $d$ is invertible. Thus if $\left[\begin{array}{cc}
a&b\\
c&d\end{array}\right]$ is also invertible, then the Schur complement $a-bd^{-1}c$ is automatically invertible  and we have
\begin{equation}\label{equation:BasicInversionFormula}
\left[\begin{array}{cc}
a&b\\
c&d\end{array}\right]^{-1}
=\left[\begin{array}{cc}
0&0\\
0&d^{-1}
\end{array}\right]+\left[\begin{array}{r}
1\\
-d^{-1}c
\end{array}\right](a-bd^{-1}c)^{-1}\left[\begin{array}{rr}
1&-bd^{-1}
\end{array}\right].
\end{equation}
Furthermore, writing
$\left[\begin{array}{cc}
a&b\\
c&d
\end{array}\right]^{-1}=\left[\begin{array}{cc}p&q\\
r&s\end{array}\right]$ with $p$, $q$, $r$, $s$ the same dimensions as $a$, $b$, $c$, $d$, respectively,
we have
\begin{equation}\label{equation:BasicInversionFormulaBis}
\left[\begin{array}{cc}
a&b\\
c&d
\end{array}\right]^{-1}=\left[\begin{array}{cc}
0&0\\
0&d^{-1}\end{array}\right]+\left[\begin{array}{c}
p\\
r\end{array}\right]p^{-1}\left[\begin{array}{cc}
p&q\end{array}\right].
\end{equation}
\subsection{Relations among the objects associated to $R$}
 We have the relation
\begin{equation}\label{equation:BasicGidNought}
\bigvee_{i=1}^N|(\{UV\}-z\Ibold_N)^{-1}(i,i)-m|\leq \bigvee_{i=1}^N\normnc{G_i-M}
\end{equation}
because the resolvent of the anticommutator  $\{UV\}$
appears as the $N$-by-$N$ block in the upper left corner of  the generalized resolvent $R$
and by definition $m=M(1,1)$. 
Let
\begin{eqnarray}\label{equation:Littlerdef}
r&=&\left[\begin{array}{ccc}
(\{UV\}-z\Ibold_N)^{-1}&0&0\\
0&0&0\\
0&0&0
\end{array}\right]\in \Mat_{3N},
\end{eqnarray}
which is just the resolvent of $\{UV\}$ bordered by some zeros.
Let $\Lambda^0\in \Mat_3$ be the constant matrix defined on line \eqref{equation:LambdaNoughtDef}.
Then we have
\begin{equation}\label{equation:BasicRidBis}
R+\Lambda^0\otimes \Ibold_N=WrW^*,\;\;\frac{dR}{dz}=Wr^2W\;\;\mbox{and}\;\;\frac{\Im R}{\Im z}=Wrr^*W=Wr^*rW
\end{equation}
as one can verify straightforwardly starting from \eqref{equation:BasicRid}. 
We  have a key {\em a priori} bound
\begin{equation}\label{equation:BasicGid}
\bigvee_{i=1}^N\normnc{G_i-M}\leq \frac{2^7(\normnc{U}\vee\normnc{V}\vee 1)^2}{\Im z}
\end{equation}
following from \eqref{equation:StupidXW},
\eqref{equation:Littlerdef},
\eqref{equation:BasicRidBis}, the definition of $G_i$
 and the bounds \eqref{equation:AnticommBounds} from Proposition \ref{Proposition:ACNondegeneracy}. The matrix $G_i$ is automatically invertible and  we have
 \begin{equation}\label{equation:KeyIdentityAC1}
-Q_i=G_i^{-1}+\Lambda+\Phi(\widehat{G}_i)
\end{equation}
by the matrix identity \eqref{equation:BasicInversionFormula}.
We furthermore have
$$
R=\hat{\ebold}_i^*R_i\hat{\ebold}_i+R\ebold_i^*G_i^{-1}\ebold_iR
$$
by the matrix identity \eqref{equation:BasicInversionFormulaBis}
and hence we have a key bound
\begin{equation}
\label{equation:KeyIdentityAC2bis}
N\normnc{G-\widehat{G}_i}\leq \normnc{G_i^{-1}}\normnc{R\ebold_i^*}_2\normnc{\ebold_iR}_2
\end{equation}
by the matrix H\"{o}lder inequality.

Here is the main result of this section.
Notably, it is a deterministic statement.

\begin{Theorem}\label{Theorem:Gizmo}
Let $\tau\geq 8$ and $\theta\geq 1$ be absolute constants.
Assume $\normnc{U}\vee \normnc{V}\leq 4$.
Consider the rectangle
\begin{equation}\label{equation:Rectangle}
\RRR=\left\{z\in \hhh\bigg\vert |\Re z|\leq 8\;\;\mbox{and}\;\; \frac{1}{N}\leq \Im z\leq \tau\right\}
\end{equation}
and let $K=\sup_{z\in \RRR}\Kfrak(z)<\infty$. (We write $\Kfrak(z)$ here to show $z$-dependence.)
Let $h$ be as defined on line \eqref{equation:hdef}.
Consider also the closed (possibly empty) set
$$
\XXX=\left\{z\in \RRR\bigg\vert  \frac{ 4 c_{\ref{Proposition:ACNondegeneracy}}^2\theta^2K^2}{N}\leq h^2\Im z\right\}
$$
where $c_{\ref{Proposition:ACNondegeneracy}}$ is the constant from Proposition \ref{Proposition:ACNondegeneracy}.
Then we have
\begin{equation}\label{equation:McGuffin}
z\in \XXX\Rightarrow \bigvee_{i=1}^N \normnc{G_i-M}\leq \frac{\theta K}{\sqrt{Nh\Im z}}\end{equation}
provided that $\tau$ is sufficiently large and $\theta$  is sufficiently large depending on $\tau$.
\end{Theorem}
\noindent We complete the proof in \S\ref{subsection:ProofOfGizmo} below.
We will prove the theorem by applying successively Propositions \ref{Proposition:Frak} and \ref{Proposition:Uroboric}.
In \S\ref{section:ProofOfMainResult}
we will construct the random variable $\Kbold$ figuring in Theorem \ref{Theorem:MainResult} by suitably approximating the quantity $\theta K$ from above.

\begin{Proposition}\label{Proposition:SelfConsistentAC}
For $i=1,\dots,N$ and $z\in\hhh$ we have
\begin{eqnarray}
\label{equation:SelfConsistentAC1}
\normnc{G_i^{-1}+\Lambda+\Phi(\widehat{G}_i)}&\leq &
4\Kfrak\normnc{W}\sqrt{\frac{(\Im z)_\bullet}{N\Im z}}\normnc{\widehat{G}_i}_\bullet^{1/2}\;\;
\mbox{and}\\
\label{equation:SelfConsistentAC2}
\normnc{G-\widehat{G}_i}&\leq & 16\normnc{W}^2\frac{(\Im z)_\bullet}{N\Im z}\normnc{G_i}_\bullet\normnc{G_i^{-1}}.
\end{eqnarray}
\end{Proposition}
 \proof 
By \eqref{equation:BasicRidBis}, \eqref{equation:BasicGid} and the matrix H\"{o}lder inequality,  we have
$$\trace \frac{\Im G_i}{\Im z}=\trace \ebold_i\frac{\Im R}{\Im z}\ebold_i^*=\normnc{\ebold_i Wr}_2^2
\geq \frac{\normnc{\ebold_iW rW^*}_2^2}{\normnc{W^*}^2}= \frac{\normnc{\ebold_i (R+\Lambda^0\otimes \Ibold_N)}_2^2}{\normnc{W}^2}$$
and similarly
$$\trace \frac{\Im G_i}{\Im z}\geq \frac{\normnc{ (R+\Lambda^0\otimes \Ibold_N)\ebold_i^*}_2^2}{\normnc{W}^2}.$$
It follows that
$$\sqrt{2}+\normnc{W}\sqrt{\trace \frac{\Im G_i}{\Im z}}
\geq\normnc{\ebold_iR}_2\vee\normnc{R\ebold_i^*}_2.$$
It follows in turn that
\begin{eqnarray}
\label{equation:Selfy1}
16\normnc{W}^2\frac{(\Im z)_\bullet}{\Im z}\normnc{G_i}_\bullet\geq 4+2\normnc{W}^2\trace\frac{\Im G_i}{\Im z}
&\geq&\normnc{\ebold_iR}_2^2\vee\normnc{R\ebold_i^*}_2^2,\\
\nonumber
16\normnc{W}^2\frac{(\Im z)_\bullet}{\Im z}\normnc{G}_\bullet\geq 4+2\normnc{W}^2\trace \frac{\Im G}{\Im z}
&\geq&\frac{\normnc{R}_2^2}{N},\;\mbox{similarly}\\
\nonumber
16\normnc{W}^2\frac{(\Im z)_\bullet}{\Im z}\normnc{\widehat{G}_i}_\bullet\geq 4+2\normnc{W}^2\trace \frac{\Im \widehat{G}_i}{\Im z}
&\geq&\frac{\normnc{R_i}_2^2}{N}\;\;\mbox{and hence}\\
\label{equation:Selfy2}
4\normnc{W}\sqrt{\frac{(\Im z)_\bullet}{N\Im z}}\normnc{\widehat{G}_i}_\bullet^{1/2}&\geq&\frac{1}{\sqrt{N}}
\left(1\vee \frac{\normnc{R_i}_2}{\sqrt{N}}\right).
\end{eqnarray}
Statements \eqref{equation:KeyIdentityAC1} and \eqref{equation:Selfy2} along with the definition of $\Kfrak$ prove \eqref{equation:SelfConsistentAC1}.
Statements \eqref{equation:KeyIdentityAC2bis} and \eqref{equation:Selfy1} prove \eqref{equation:SelfConsistentAC2}.
\qed

\begin{Proposition}\label{Proposition:FrakApplication}
For every $z\in \hhh$ we have
\begin{equation}\label{equation:BigBump}
\bigvee_{i=1}^N \normnc{G_i-M}\leq \frac{\sqrt{h}}{c_{\ref{Proposition:ACNondegeneracy}}}\Rightarrow
\bigvee_{i=1}^N \normnc{G_i-M}\leq \frac{C(8+|z|)^5\normnc{W}\Kfrak}{\sqrt{Nh\Im z}}
\end{equation}
where  $C$ is an absolute constant.
\end{Proposition}
\proof Proposition \ref{Proposition:Frak} specialized to the present setup is the assertion that
\begin{eqnarray*}
&&\bigvee_{i=1}^N \normnc{G_i-M}\leq \frac{1}{8\normnc{\kappa}_\bullet\normnc{\Phi}_\bullet}\\
\nonumber&\Rightarrow &\bigvee_{i=1}^N \normnc{G_i-M}\leq 2^{14}(1+\normnc{M})^7(\normnc{\Phi}_\bullet\vee\normnc{\Lambda}_\bullet)^4\normnc{\kappa}_\bullet\Efrak
\end{eqnarray*}
where  the quantity $\Efrak$ satisfies
$$\Efrak\leq 4\Kfrak\normnc{W}\sqrt{\frac{(\Im z)_\bullet}{N\Im z}}$$
by Proposition \ref{Proposition:SelfConsistentAC} and the definition of $\Kfrak$.
We obtain \eqref{equation:BigBump} after simplifying by means of Proposition \ref{Proposition:ACNondegeneracy}.
\qed


\subsection{Proof of Theorem \ref{Theorem:Gizmo}}
\label{subsection:ProofOfGizmo}
On the set $\XXX$ we consider the three continuous functions
$$f_1=\bigvee_{i=1}^N \normnc{G_i-M}\bigg\vert_\XXX,\;\;f_2=\frac{\sqrt{h}}{c_{\ref{Proposition:ACNondegeneracy}}}\bigg\vert_\XXX\;\;
\mbox{and}\;\;
f_3= \frac{\theta K}{\sqrt{Nh\Im z}}\bigg\vert_\XXX.
$$
The rest of the proof is a matter of checking hypotheses in  Proposition \ref{Proposition:Uroboric}.

\subsubsection{$\XXX$ if not empty is connected}
Let 
$$\rho=\frac{4c_{\ref{Proposition:ACNondegeneracy}}^2\theta^2K^2}{N}.$$
Since $\tau\geq 1$ and $\rho\geq \frac{1}{N}$ by assumption, and $h\equiv 1$ on the imaginary axis,
the set $\XXX$ is not empty if and only if $\rho\leq \tau$.
It follows that $\XXX$ is nonempty if and only if it contains the horizontal line segment
$$\ii\tau+[-8,8]=\{x+\ii\tau\mid -8\leq x\leq 8\}.$$
Furthermore, since  the function $h^2\Im z$ is monotone increasing on vertical lines in $\hhh$,
each point of $\XXX$ is connected to the line segment $\ii\tau+[-8,8]$
by a vertical line segment contained in $\XXX$.  
Thus, indeed, $\XXX$ if not empty is connected.

\subsubsection{
Checking  hypothesis \eqref{equation:Uroboric1}}
Using \eqref{equation:StupidXW},  \eqref{equation:BasicGid} and our hypothesis $\normnc{U}\vee \normnc{V}\leq 4$
to justify the first inequality below, we choose any $\tau\geq 8$ large enough to make
the statement
$$
\bigvee_{i=1}^N \normnc{G_i-M}\bigg\vert_{z=\ii\tau}\leq \frac{2^{11}}{\Im z}\bigg\vert_{z=\ii\tau}
=\frac{2^{11}}{\tau}<\frac{1}{c_{\ref{Proposition:ACNondegeneracy}}}=
\frac{\sqrt{h}}{c_{\ref{Proposition:ACNondegeneracy}}}\bigg\vert_{z=\ii\tau}$$
hold.
With $\tau$ thus fixed, hypothesis \eqref{equation:Uroboric1}  of Proposition \ref{Proposition:Uroboric}
is verified.
\subsubsection{
Checking  hypothesis \eqref{equation:Uroboric2}}
We next choose $\theta$ so that
$$\theta\geq C_{\ref{Proposition:FrakApplication}}(16+\tau)^52^5.$$
Then by Proposition \ref{Proposition:FrakApplication}, \eqref{equation:StupidXW}
and our hypothesis that $\normnc{U}\vee \normnc{V}\leq 4$, we have
$$
z\in \XXX\;\;\mbox{and}\;\;\bigvee_{i=1}^N \normnc{G_i-M}\leq \frac{\sqrt{h}}{c_{\ref{Proposition:ACNondegeneracy}}}\Rightarrow
\bigvee_{i=1}^N \normnc{G_i-M}\leq \frac{\theta K}{\sqrt{Nh\Im z}}.
$$
With $\theta$ thus fixed, hypothesis \eqref{equation:Uroboric2}  of Proposition \ref{Proposition:Uroboric} is verified.
\subsubsection{
Checking  hypothesis \eqref{equation:Uroboric3}}
Finally we have
$$
z\in \XXX\Rightarrow \frac{\theta K}{\sqrt{Nh\Im z}}\leq \frac{\sqrt{h}}{2c_{\ref{Proposition:ACNondegeneracy}}}<\frac{\sqrt{h}}{c_{\ref{Proposition:ACNondegeneracy}}}
$$
by the very definition of $\XXX$. 
Thus hypothesis \eqref{equation:Uroboric3}  of Proposition \ref{Proposition:Uroboric} is verified. 
The conclusion \eqref{equation:Uroboric4} of Proposition \ref{Proposition:Uroboric} is then
the same as the conclusion \eqref{equation:McGuffin} of Theorem \ref{Theorem:Gizmo}. \qed

The following technical assertion will be needed in the next section.
We write $\Kfrak_i(z)$ to show $z$-dependence.
\begin{Proposition}\label{Proposition:Lipschitz}
Assume that $\normnc{U}\vee \normnc{V}\leq 4$. 
For $i=1,\dots,N$ 
and $z_1,z_2\in \hhh$ such that $(\Im z_1)\wedge(\Im z_2)\geq \frac{1}{N}$
we have
$$\frac{\left|\Kfrak_i(z_1)-\Kfrak_i(z_2)\right|}{|z_1-z_2|}\leq cN^{7/2},$$
where $c$ is an absolute constant. 
\end{Proposition}
\proof  
The proof is just an extremely ugly computation based on \eqref{equation:BasicRidBis}.
Let $$C=1+\normnc{\Phi}+\normnc{X}+\normnc{X}^2+3\normnc{W}^2.$$
Temporarily (only in this proof) we write $R(z)$, $R_i(z)$ and $Q_i(z)$ when necessary to show $z$-dependence.
We evidently have 
$$\normnc{R}\leq 2+\normnc{W}^2N\leq 3\normnc{W}^2N\leq CN\;\;
\mbox{and similarly}\;\;
\normnc{R_i}\leq CN.$$
Consequently we have
$$\normnc{Q_i}\leq (\normnc{X}^2+\normnc{X}+\normnc{\Phi})\normnc{R_i}\leq C^2N.$$
We may assume that $z_1\neq z_2$. We have
$$\frac{\normnc{R(z_1)-R(z_2)}}{|z_1-z_2|}\leq N^2\normnc{W}^2\leq CN^2,\;\;
\mbox{similarly}\;\;
\frac{\normnc{R_i(z_1)-R_i(z_2)}}{|z_1-z_2|}\leq CN^2,$$
hence
$$\frac{\normnc{R_i(z_1)-R_i(z_2)}_2}{|z_1-z_2|}\leq  C N^{5/2}$$
and furthermore
$$\frac{\normnc{Q_i(z_1)-Q_i(z_2)}}{|z_1-z_2|}\leq 
(\normnc{X}^2+\normnc{\Phi})CN^2\leq C^2N^2.$$
Now consider the functions
$$f=\sqrt{N}\vee \normnc{R_i}_2\bigg\vert_{\{\Im z\geq \frac{1}{N}\}}\;\;\mbox{and}\;\;\;\;
g=N\normnc{Q_i}\bigg\vert_{\{\Im z\geq \frac{1}{N}\}}.$$
To finish the proof we must estimate the Lipschitz constant of $1\vee \frac{g}{f}$
and thus need only estimate that of $\frac{g}{f}$.
We have thus far determined that $f$ is (upper) bounded by $CN$
and has Lipschitz constant bounded by $CN^{5/2}$;
also by definition $f$ is lower bounded by $\sqrt{N}$.
Furthermore we have determined that $g$ is bounded by $C^2N^2$ and has Lipschitz constant bounded by $C^2N^3$.
Using the identity
$$\frac{g(z_1)}{f(z_1)}-\frac{g(z_2)}{f(z_2)}=(g(z_1)-g(z_2))\frac{1}{f(z_1)}+g(z_2)\frac{f(z_2)-f(z_1)}{f(z_1)f(z_2)}
$$
we deduce that
$$
\frac{\left|\frac{g(z_1)}{f(z_1)}-\frac{g(z_2)}{f(z_2)}\right|}{|z_1-z_2|}
\leq C^2N^3/\sqrt{N}+(C^2N^2)(CN^{5/2})/N\leq2C^3N^{7/2},
$$
which finishes the proof.
\qed

\section{Proof of Theorem \ref{Theorem:MainResult}}\label{section:ProofOfMainResult}
In this section we work simultaneously in the settings of Theorem \ref{Theorem:MainResult}
and Theorem \ref{Theorem:Gizmo}.  We fix once and for all absolute constants $\tau\geq 8$ and $\theta\geq 1$ 
so that the conclusion \eqref{equation:McGuffin} of Theorem \ref{Theorem:Gizmo} holds.

\subsection{Construction of $\Kbold$}
 Let $\RRR$ be the rectangle \eqref{equation:Rectangle}.
By Proposition \ref{Proposition:Lipschitz} we know that conditioned on $\normnc{U}\vee \normnc{V}\leq 4$
the quantity $\Kfrak_i(z)$ depends Lipschitz-continuously
on $z\in \RRR$  with Lipschitz constant bounded by $cN^{7/2}$. Thus for suitable absolute constants
$\beta_1$ and $\beta_3$ and a suitable net $\RRR_0\subset \RRR$ of at most $\beta_3N^{\beta_1-1}$ points we have
\begin{equation}\label{equation:RectangleBis}
 2\bigvee_{i=1}^N \bigvee_{z_0\in \RRR_0}\Kfrak_i(z_0)
 \geq \sup_{z\in \RRR}\Kfrak(z)
\end{equation}
conditioned on $\normnc{U}\vee\normnc{V}\leq 4$.
We define $\Kbold$ to equal the left side of \eqref{equation:RectangleBis} multiplied by $\theta$.
By Theorem \ref{Theorem:Gizmo} and the bound \eqref{equation:BasicGidNought}
the random variable $\Kbold\geq 1$ automatically has property
\eqref{equation:FirstKboldProperty}. It remains only to prove that $\Kbold$ has property
\eqref{equation:SecondKboldProperty}.
The latter task is just a matter of revisiting the topic of \cite[Appendix B]{EYY},
namely large deviations for quadratic forms in independent variables satisfying exponential tail bounds.
However, because we have to make a few adjustments
to handle the special features of our anticommutator setup,
we will handle the details a bit differently than in the cited reference.

\subsection{Remark} In the proof of the local semicircle law \cite[Thm. 3.1]{EYY} the Lipschitz continuity of the various functions in play 
is frequently invoked while marching toward the real axis.
It might have seemed we were trying to avoid such considerations here
by using Proposition \ref{Proposition:Uroboric}. Certainly we have avoided their use in a dynamical  way.
But ultimately our reworking of the method of \cite{EYY} has merely displaced the use of Lipschitz continuity to the phase of the argument presented immediately above in which we  construct the random variable $\Kbold$.\\

We begin the proof that $\Kbold$ has property \eqref{equation:SecondKboldProperty} by
recalling the simple relationship between moment bounds of the form \eqref{equation:Wig4}
and exponentially light tails.

\begin{Proposition}\label{Proposition:EasyDeavy}
Fix constants $\alpha>0$ and $C\geq 1$.
Let $Z$ be a nonnegative random variable.
(i) If 
$\sup_{p\in [2,\infty)}p^{-\alpha}\norm{Z}_p\leq 1$
then 
$\Pr(Z>t^\alpha)\leq  \exp\left(\alpha\left(2-\frac{t}{e}\right)\right)$.
(ii) If $\Pr(Z>t^\alpha)\leq Ce^{-t}$ for $t>0$
then $\sup_{p\in [2,\infty)}p^{-\alpha}\norm{Z}_p
\leq C^{1/2}\left(\alpha+\frac{1}{2}\right)^\alpha$.
\end{Proposition}
\proof (i) In the Markov bound 
$\Pr(Z>t^\alpha)\leq \frac{\norm{Z}_p^p}{t^{\alpha p}}\leq \frac{p^{\alpha p}}{t^{\alpha p}}$
we substitute $p=t/e$ if $t/e\geq 2$ and simplify. 
(ii) 
For the $\Gamma$-function $\Gamma(s)=\int_0^\infty x^{s-1}e^{-x}\,dx$ one has a functional equation
$s\Gamma(s)=\Gamma(s+1)$, a bound $\Gamma(s)\leq 1$ for $1\leq s\leq 2$ and (hence) an elementary inequality
$\Gamma(1+s)\leq (1+s)^s$ for $s\geq 0$.
For $p\geq 1$ we then have
$$
\Ebold Z^p
=\alpha p\int_0^\infty \Pr(Z>t^\alpha)\,t^{\alpha p-1}\,dt
\leq \alpha pC\int_0^\infty e^{-t}\,t^{\alpha p-1}\,dt\;\leq C(1+p\alpha )^{p\alpha }
$$
and thus $p^{-\alpha}\norm{Z}_p\leq C^{1/p}\left(\alpha +\frac{1}{p}\right)^\alpha$
for $p\geq 1$. 
\qed

We next recall a classical result.
Let $\Theta(s)=\frac{2^{s/2}}{\sqrt{\pi}}\Gamma\left(\frac{s+1}{2}\right)$
for $s\geq 0$.

\begin{Theorem}[Whittle \cite{Whittle}]\label{Theorem:Whittle}
 Let
$Y_1,\dots,Y_n$ be independent real random variables of mean zero.
Fix $p\in [2,\infty)$.
Let $v\in \RR^n$ be a real vector.
Let $B\in \Mat_n$ be a matrix with real entries.
If $\bigvee_{i=1}^n \norm{Y_i}_p<\infty$, then
\begin{equation}\label{equation:WhittleLinear}
\norm{\sum_{i=1}^n v(i)Y_i}_p
\leq 2\Theta(p)^{\frac{1}{p}}\left(\sum_{i=1}^n v(i)^2\norm{Y_i}_p^2\right)^{1/2}.
\end{equation}
Furthermore, if $\bigvee_{i=1}^n \norm{Y_i}_{2p}<\infty$, then
\begin{eqnarray}\label{equation:WhittleQuadratic}
&&\norm{\sum_{i,j=1}^N B(i,j)(Y_iY_j-\Ebold[Y_iY_j])}_p\\
\nonumber&\leq& 2^3\Theta(p)^{\frac{1}{p}}\Theta(2p)^{\frac{1}{2p}}\left(\sum_{i,j=1}^N B(i,j)^2\norm{Y_i}_{2p}^2\norm{Y_j}^2_{2p}\right)^{1/2}.
\end{eqnarray}
\end{Theorem}
\noindent We hasten to point out that one has an elementary bound
\begin{equation}\label{equation:SimpleGammaBound}
\sup_{s\geq 2}\frac{\Theta(s)^{\frac{1}{s}}}{\sqrt{s}}\leq 1.
\end{equation}
Thus the estimates \eqref{equation:WhittleLinear}
and \eqref{equation:WhittleQuadratic} can be simplified nicely. 

From Proposition \ref{Proposition:EasyDeavy}
and Theorem \ref{Theorem:Whittle} we then get the following tail-bound.
\begin{Proposition}\label{Proposition:QuadEasyDeavy}
Fix constants $\gamma_0>0$ and $\gamma_1\geq 1$. Fix a positive integer $k$.
Let $Y_0,\dots,Y_{2N}\in \Mat_k$ be random of mean zero.
Assume that the family 
$$\{\sigma(Y_0)\}\cup\left\{\sigma(Y_i,Y_{i+N})\right\}_{i=1}^N$$ of $\sigma$-fields is independent.
Assume that
$$\sup_{p\geq 2}p^{-\gamma_0}\bigvee_{i=0}^{2N}\norm{\normnc{Y_i}}_p\leq \sqrt{\frac{\gamma_1}{N}}.$$
Let 
$$Y=\left[\begin{array}{ccc}Y_1&\dots& Y_N\end{array}\right]\in \Mat_{k\times kN}\;\;\mbox{and}\;\;
\widehat{Y}=\left[\begin{array}{ccc}Y_{N+1}&\dots&Y_{2N}\end{array}\right]\in \Mat_{k\times kN}.$$
Let $B\in \Mat_{kN}$ be any constant matrix.
Then for every $t>0$ we have
$$
\Pr\left(\frac{\normnc{YB\widehat{Y}^*-Y_0-\Ebold(YB\widehat{Y}^*)}}
{\frac{\gamma_1}{\sqrt{N}}\left(1\vee \frac{\normnc{B}_2}{\sqrt{N}}\right)}>t^{2\gamma_0+1} \right)
\leq \gamma_2e^{-\gamma_3t}
$$
for constants $\gamma_2\geq 1$ and $\gamma_3>0$ depending only on $\gamma_0$ and $k$.
\end{Proposition}
\proof  After  replacing $Y_0$ by $Y_0/\gamma_1$
and for $i>0$ replacing $Y_i$ by $Y_i/\sqrt{\gamma_1}$,
we may assume without loss of generality that $\gamma_1=1$.
By Proposition \ref{Proposition:EasyDeavy} we already have a tail bound
$\Pr(\sqrt{N}\normnc{Y_0}>t^{\gamma_0})\leq Ce^{-ct}$.
Thus we may assume that $Y_0=0$.
In turn by Proposition \ref{Proposition:EasyDeavy} it is enough to prove
that
\begin{equation}\label{equation:DesiredBound}
\sup_{p\geq 2}p^{-(1+2\gamma_0)}\norm{\normnc{YB\widehat{Y}^*-\Ebold(YB\widehat{Y}^*)}}_p\leq \frac{\gamma_4
\normnc{B}_2}{N}\end{equation}
where $\gamma_4\geq 1$ is a constant depending only on $\gamma_0$ and $k$. 
Without loss of generality we may assume that
$B$ has real entries and that the random matrices $Y_i$ have real entries.
We may then in turn obviously assume that $k=1$.
By \eqref{equation:WhittleLinear} we may
assume that every diagonal entry of $B$ vanishes.
We may also obviously assume that $N\geq 2$.
Now let $I\subset \{1,\dots,N\}$ be any subset of cardinality
$\lfloor\frac{N}{2}\rfloor$ and let $I^c$ denote the complement of $I$.
Let $$B_I(i,j)=B(i,j)\one_{i\in I}\one_{j\in I^c},$$
thus defining a matrix $B_I\in \Mat_N$
supported on the set
$$I\times I^c\subset \{1,\dots,N\}^2.$$
Let 
$$\widetilde{Y}_I(i)=\left\{\begin{array}{rl}
Y_i&\mbox{if $i\in I$,}\\
Y_{i+N}&\mbox{if $i\in I^c$.}
\end{array}\right.$$
Note that the entries of $\tilde{Y}_I$ are independent.
Note also that 
$$YB_I\widehat{Y}^*=\widetilde{Y}_IB_I\widetilde{Y}_I^*.$$
Thus we have
\begin{equation}\label{equation:DesiredBoundBis}
\sup_{p\geq 2}p^{-(1+2\gamma_0)}\norm{\normnc{YB_I\widehat{Y}^*-\Ebold(YB_I\widehat{Y}^*)}}_p\leq \frac{\gamma_4}{4}\frac{\normnc{B_I}_2}{N}\leq \frac{\gamma_4}{4}\frac{\normnc{B}_2}{N}
\end{equation}
by Theorem \ref{Theorem:Whittle}
and the upper bound \eqref{equation:SimpleGammaBound}. 
Now the average of $B_I$ over $I$ equals $qB$ for some constant $q\geq \frac{1}{4}$.
Thus,  averaging over $I$ on the left side of \eqref{equation:DesiredBoundBis} and using
Jensen's inequality, we obtain \eqref{equation:DesiredBound}.
 \qed

 \subsection{End of the proof}
 The summary of properties of the random matrix $X$
 in \S\ref{subsection:PhiRecog} and Proposition \ref{Proposition:QuadEasyDeavy} together  provide us with constants $\beta_2\geq 1$
 and $\beta_4>0$ depending only on $\alpha_0$ and $\alpha_1$
 such that for $i=1,\dots,N$ and any $z_0\in \hhh$ we have a uniform conditional tail bound
 $$\Pr\left(2\theta\Kfrak_i(z_0)>t^{\frac{1}{2\alpha_0+1}}\bigg\vert \hat{\ebold}_iX\hat{\ebold}_i^*\right)\geq \beta_4\exp(-\beta_2t)\;\;\mbox{a.s..}
 $$
The latter combined with the evident union bound over $\beta_3N^{\beta_1}$ events yields \eqref{equation:SecondKboldProperty}
with $\beta_0=\beta_3\beta_4$.
The proof of Theorem \ref{Theorem:MainResult} is complete. \qed

\section{Appendix: The deterministic local semicircle law}\label{section:Appendix}
We state and sketch the proof of the semicircular analogue of Theorem \ref{Theorem:Gizmo}.
For the proof we will use Propositions \ref{Proposition:Uroboric},  \ref{Proposition:SemicircleReview}
and  \ref{Proposition:Frak}, none of which have anything to do specifically with anticommutators.
\subsection{Setup for the result}
\subsubsection{Basic data}
Fix a hermitian matrix  $X\in \Mat_N$.
\subsubsection{Specialized matrix notation}
Let $e_i$ denote the $i^{th}$ row of $\Ibold_N$
and let $\hat{e}_i$ denote the result of deleting the $i^{th}$ row of $\Ibold_N$.
\subsubsection{Functions of $z$}
For $z\in \hhh$ and $i=1,\dots,N$ let 
\begin{eqnarray*}
R&=&(X-z\Ibold_N)^{-1}\in \Mat_N,\;\;G_i\;=\;e_iRe^*_i=R(i,i)\in \hhh,\\
R_i&=&(\hat{e}_iX\hat{e}_i^*-z\Ibold_{N-1})^{-1}\in \Mat_{N-1},\;\;
\widehat{G}_i\;=\;\frac{1}{N}\trace R_i\in \hhh,\\
Q_i&=&
e_iX\hat{e}_i^*R_i\hat{e}_iXe_i^*-X(i,i)-\widehat{G}_i\in \CC,\\
m&=&\frac{1}{2\pi}\int_{-2}^2 \frac{\sqrt{4-t^2}\,dt}{t-z}\in \hhh\;\;(\mbox{equivalently: $z=-m^{-1}-m$}),\\
\label{equation:hdefSC}
h&=&1\wedge |z+2|\wedge |z-2|>0,\\
\Kfrak&=&
 1\vee\bigvee_{i=1}^N
\frac{|Q_i|}{\frac{1}{\sqrt{N}}\left(1\vee \frac{\normnc{R_i}_2}{\sqrt{N}}\right)}\in [1,\infty).
\end{eqnarray*}
All these objects depend on $z$ but the notation does not show it.
Note that the $z$-dependence is continuous.

Here then is the {\em deterministic local semicircle law}.
\begin{Theorem}\label{Theorem:GizmoSC} Notation and assumptions are as above.
Let $\tau=20$ and $\theta=2^{100}$.
Consider the rectangle
\begin{equation}\label{equation:RectangleSC}
\RRR=\left\{z\in \hhh\bigg\vert |\Re z|\leq 4\;\;\mbox{and}\;\; \frac{1}{N}\leq \Im z\leq \tau\right\}
\end{equation}
and let  $K=\sup_{z\in \RRR}\Kfrak(z)<\infty$.
(We write $\Kfrak(z)$ to show $z$-dependence.)
Consider also the closed (possibly empty) set
$$
\XXX=\left\{z\in \RRR\bigg\vert\frac{2^8\theta^2 K^2}{N}\leq h^2\Im z\right\}.
$$
Then we have
\begin{equation}\label{equation:McGuffinSC}
z\in \XXX\Rightarrow \bigvee_{i=1}^N |G_i-m|\leq \frac{\theta K}{\sqrt{Nh\Im z}}.
\end{equation}
\end{Theorem}
\noindent We break the proof into several stages.
\subsection{Application of Propositions \ref{Proposition:SemicircleReview} and \ref{Proposition:Frak}}
Let
$$
\Efrak=\bigvee_{i=1}^N
\frac{|G_i^{-1}+z+\widehat{G}_i|}{|\widehat{G}_i|_\bullet^{1/2}}\vee\bigvee_{i=1}^N
\sqrt{\frac{|\widehat{G}_i-G|}{|G_i|_\bullet|G_i^{-1}|}.}
$$
To use Proposition \ref{Proposition:Frak} the main thing is to bound $\Efrak$.  In any case we have
$$
\frac{\normnc{R_i}_2^2}{N}=\frac{\Im \widehat{G}_i}{\Im z}\;\;\mbox{and}\;\;
\normnc{e_iR}_2\vee \normnc{Re_i^*}_2=\frac{\Im G_i}{\Im z}.
$$
Both statements are merely rewrites of \eqref{equation:MagicResolvent} above.
Thus, since 
$$Q_i=-G_i^{-1}-z-\widehat{G}_i$$ by \eqref{equation:BasicInversionFormula}, we have
\begin{eqnarray*}
|G_i^{-1}+z+\widehat{G}_i|=|Q_i|&\leq &\frac{\Kfrak}{\sqrt{N}}
\left(1\vee \sqrt{\frac{\Im G_i}{\Im z}}\right)\leq \Kfrak\sqrt{\frac{(\Im z)_\bullet}{N\Im z}}|\widehat{G}_i|_{\bullet}^{1/2},\\
|G-\widehat{G}_i|&\leq & \frac{1}{N\Im z}\leq \frac{(\Im z)_\bullet }{N\Im z}|G_i|_\bullet |G_i|^{-1}
\end{eqnarray*}
and hence 
$$
\Efrak\leq \sqrt{\frac{(\Im z)_\bullet}{N\Im z}}\Kfrak.
$$
Let $\kappa=(m^{-1}-m)^{-1}$. By Proposition \ref{Proposition:SemicircleReview} the quadruple $(z,m,1,\kappa)$
is a nondegenerate solution of the Schwinger-Dyson equation defined over $\CC$.
Thus we have
$$
\bigvee_{i=1}^N |G_i-m|\leq \frac{1}{8|\kappa|_\bullet}\Rightarrow
\bigvee_{i=1}^N |G_i-m|\leq   \frac{2^{21}|z|_\bullet^4\sqrt{(\Im z)_\bullet}\Kfrak}{\sqrt{N\Im z}}|\kappa|_\bullet
$$
by substituting into Proposition \ref{Proposition:Frak}.
In turn we obtain the statement
\begin{equation}\label{equation:AlmostThere}
\bigvee_{i=1}^N |G_i-m|\leq \frac{\sqrt{h}}{8}\Rightarrow
\bigvee_{i=1}^N |G_i-m|\leq   \frac{2^{21}|z|_\bullet^5\Kfrak}{\sqrt{Nh\Im z}}
\end{equation}
after using the bound $\normnc{\kappa}_\bullet \leq\frac{1}{\sqrt{h}}$
noted in Proposition \ref{Proposition:SemicircleReview} and simplifying slightly. 

\subsection{Checking hypotheses in Proposition \ref{Proposition:Uroboric}}

On the set $\XXX$ we consider the three continuous functions
$$f_1=\bigvee_{i=1}^N |G_i-M|\bigg\vert_\XXX,\;\;f_2=\frac{\sqrt{h}}{8}\bigg\vert_\XXX\;\;
\mbox{and}\;\;
f_3= \frac{\theta K}{\sqrt{Nh\Im z}}\bigg\vert_\XXX.
$$
We now have only to check hypotheses in  Proposition \ref{Proposition:Uroboric}
in order to finish the proof of Theorem \ref{Theorem:GizmoSC}.

\subsubsection{$\XXX$ is connected if nonempty}
Let 
$$\rho=2^8\theta^2K^2/N.$$
 Note that $\rho\geq 1/N$ and $\tau\geq 1$.
 In terms of the parameters $\rho$ and $\tau$ we have
$$\XXX=\{z\in \hhh\mid |\Re z|\leq 4,\;\Im z\leq \tau\;\mbox{and}\;\rho\leq h^2\Im z\}.$$
Thus a necessary and sufficient condition for nonemptiness is that $\rho\leq \tau$,
and under those equivalent conditions $\XXX$ automatically contains the line segment $\ii\tau+[-4,4]$.
Finally, since $h^2\Im z$ is monotone increasing on vertical lines in $\hhh$,
every point of $\XXX$ is connected by a vertical line segment in $\XXX$ to $\ii\tau+[-4,4]$.

\subsubsection{Checking hypothesis \eqref{equation:Uroboric1}}
We have evident bounds
$$
|m|\vee \bigvee_{i=1}^N |G_i|\leq \frac{1}{\Im z},\;\;\mbox{hence}\;\;\bigvee_{i=1}^N |G_i-m|\leq \frac{2}{\Im z},
$$
hence 
$$\frac{2}{\Im z}\bigg\vert_{z=\ii \tau}=\frac{2}{\tau}<
\frac{1}{8}=
\frac{\sqrt{h}}{8}\bigg\vert_{z=\ii \tau}$$
because $\tau$ is large enough
and hence
$$\bigvee_{i=1}^N|G_i-m|\bigg\vert_{z=\ii\tau}
<\frac{\sqrt{h}}{8}\bigg\vert_{z=\ii \tau}.
$$
Thus hypothesis \eqref{equation:Uroboric1} of Proposition \ref{Proposition:Uroboric} holds.

\subsubsection{Checking hypothesis \eqref{equation:Uroboric2}}
Because $\theta$ is large enough
we have
$$\theta\geq 2^{21}|4+\tau|_\bullet^5\geq \sup_{z\in \RRR}2^{21}|z|_\bullet^5$$
and hence
$$
\left(z\in \RRR\;\mbox{and}\;\bigvee_{i=1}^N|G_i-m|\leq\frac{\sqrt{h}}{8}\right)
\Rightarrow \bigvee_{i=1}^N |G_i-m|\leq \frac{\theta K}{\sqrt{Nh\Im z}}
$$
by \eqref{equation:AlmostThere}.
Thus hypothesis \eqref{equation:Uroboric2} of Proposition \ref{Proposition:Uroboric} holds.
\subsubsection{Checking hypothesis \eqref{equation:Uroboric3}}
We have
$$
z\in \XXX\Rightarrow \frac{\theta\Kfrak}{\sqrt{Nh\Im z}}\leq \frac{\sqrt{h}}{16}<\frac{\sqrt{h}}{8}
$$
by the very definition of $\XXX$.  Thus hypothesis \eqref{equation:Uroboric2} of Proposition \ref{Proposition:Uroboric} holds.
The conclusion \eqref{equation:Uroboric4} of Proposition  \ref{Proposition:Uroboric} and conclusion
\eqref{equation:McGuffinSC}
of Theorem \ref{Theorem:GizmoSC} are then the same. The proof of Theorem \ref{Theorem:GizmoSC} is complete. \qed

\subsection{Remark}
By studying the generalized resolvent
$$\left[\begin{array}{cc}
-z\Ibold_p&X\\
X^*&-\Ibold_q
\end{array}\right]^{-1}\;\;\;(X\in \Mat_{p\times q})$$
one can obtain a similar deterministic local Marcenko-Pastur law.



\textbf{Acknowledgments} I thank J. Yin for patient explanations given to me at Oberwolfach in May 2011.
I also learned a lot from participants of the IMA workshop in June 2012.
I thank L. Erd\"{o}s and O. Zeitouni for valuable advice.

\end{document}